\documentclass{article}
\usepackage{amssymb}
\usepackage{amsfonts}
\usepackage{amsmath}

\input{tcilatex}
\begin{document}

\title{On *-fusion frames for Hilbert C*-modules }
\author{Nadia Assila \and Samir Kabbaj \and Hicham Zoubeir}
\maketitle

\begin{abstract}
Our main goal in this paper, is to generalize to Hilbert C*-modules the
concept of fusion frames. Indeed we introduce the notion of *--fusion frames
associated to weighted sequences of orthogonally complemented submodules of
a Hilbert C*-module, and prove for such a *-fusion frames some fundamental
results. 
\begin{equation*}
\text{\textbf{2010 Mathematics Subject Classification} : 42C15, 46L05, 46L08.%
}
\end{equation*}%
\begin{equation*}
\text{\textbf{Key words}}:\text{C*-algebras, HilbertC*-module, *-fusion
frames.}
\end{equation*}
\end{abstract}

\section{Introduction}

In 1946, Gabor (\cite{GABO}) performed a new method for the signal
reconstruction from elementary signals. In 1952, Duffin and Schaeffer (\cite%
{DUFF}) developped, in the field of nonharmonic series, a similar tool and
introduced frame theory for Hilbert spaces. For more than thirty years, the
results of Duffin and Schaeffer has not received from the mathematical
community, the interest they deserve, until the publication of the work of
Higgins and Young (\cite{YOUN}) where the authors studied frames in abstract
Hilbert spaces. In 1986, the work of Daubechies, Grossmann and Meyer (\cite%
{DAUB}) gave to frame theory the momentum it lacked and allowed it to be
widely studied. From this date, several generalizations of the notion of
frame have been developed, for example the concept of atomic decompositions (%
\cite{FEIC}), the concept of continuous frames in Hilbert spaces (\cite{ALI}%
), the concept of frames in Banach spaces (\cite{CASA1}), the notion of $p$%
-frames ((\cite{ALDR}), the notion of frames associated with measurable
spaces (\cite{GABA}), (\cite{CHRI})), the concept of frames of subspaces (%
\cite{CASA2}). In 2000, Frank and Larson (\cite{FRAN}) introduced the notion
of frames in Hilbert C*-modules as a generalization of frames in Hilbert
spaces. Since this innovative work, several generalizations of Hilbertian
frame theory have emerged. Let us cite for example the work of A. Khosravi
and B. Khosravi (\cite{KHOS}) who introduced the notion of g-frames in a
Hilbert C*-module, the contribution of A. Alijani and M. Dehghan (\cite{ALIJ}%
) who developed the notion of *- frame in a Hilbert C*-module, the work of
N. Bounader and S.Kabbaj (\cite{BOUN}) who introduced the notion of
*-g-frame for Hilbert C*-modules. Fusion frames as a generalisation of
frames were introduced by Casazza, Kutyniok and Li in (\cite{CASA3}). The
motivation behind fusion frame comes from the need to efficiently process
and analyze large data sets. A natural idea is to divide these data sets
into suitable smaller and parallel ones that can be processed independently.
Fusion frame systems are created to meet these needs through the weighted
coherent combination of subsystems. By this approach the fusion frame theory
provides a flexible framework that takes into consideration local frames
associated with the subsystems (\cite{CASA3}).

Our main goal in this paper, is to generalize to Hilbert C*-modules the
concept of fusion frames. Indeed we introduce the notion of *--fusion frames
associated to weighted sequences of orthogonally complemented submodules of
a Hilbert C*-module, and prove for such a *-fusion frames some fundamental
results.

The paper is structured as follows :

In section 2, we state some notations definitions and some results that are
useful for the for the proofs of the fundamental results of the paper.

In section 3, we give the definition of *-fusion frame associated to a
weighted sequence of orthogonally complemented submodules of a Hilbert
C*-module and the definition of the analysis and synthesis operators related
to a *-fusion frame. After that we prove the following properties of
*-fusion frames :

\begin{itemize}
\item The analysis operator of a *-fusion frama is a well-defined bounded
linear adjointable operator.

\item The synthesis operator of a *-fusion frame is a well-defined bounded
linear positive and self-adjoint and invertible operator. Furthermore this
operator gives rise to a reconstruction formula.

\item Every Hilbert C*-module which has a *-fusion frame is countably
generated.

\item There exists a characterisation of *-fusion frames in term of the norm
of the Hilbert C*-module.

\item Let $E$ and $F$\ be a Hilbert $\mathfrak{A}$-modules. Under some
conditions on $E$, the image of a *-fusion frame of $E$ by an orthogonally
perserving mapping from a the Hilbert $\mathfrak{A}$-module $E$ to the
Hilbert $\mathfrak{A}$-module $F$ is also a *-fusion frame of $F.$

\item Let $\left( \mathfrak{H}_{n}\right) _{n\in 
\mathbb{N}
^{\ast }}$ be a sequence of orthogonally complemented submodules of a
Hilbert $\mathfrak{A}$-module $\mathfrak{H.}$ The set of weights $\left(
\omega _{n}\right) _{n\in 
\mathbb{N}
^{\ast }}$such that\textit{\ }$\left( \left( \mathfrak{H}_{n},\omega
_{n}\right) \right) _{n\in 
\mathbb{N}
^{\ast }}$\textit{\ }is a *-fusion frame of $\mathfrak{H,}$ is a convex cone
of the $%
\mathbb{C}
$-vector space $\mathfrak{A}^{%
\mathbb{N}
^{\ast }}.$
\end{itemize}

Section 4 is devoted to the the detailed study of two examples of *-fusion
frames.

In section 5 :

\begin{itemize}
\item using a distance already introduced by Dragan S. Djordjevic in (\cite%
{DJOR}), we define the notion of the angle of two orthogonally complemented
submodules of a Hilbert C*-module;

\item using the same distance, we introduce on the set of all\ the sequences
of orthogonally complemented submodules of a Hilbert C*-module a topology
defined by an ecart;

\item relying on the constructed topolgy and on the notion of angle that we
have introduced, we prove for *-fusion frames, a perturbation results of
topological and geometric character.
\end{itemize}

For all the material on C*-algebras and Hilbert C*-modules, one can refer to
the references (\cite{DAVI}), (\cite{LANC}), (\cite{MANU}) and (\cite{MURP}).

\bigskip

\section{Preliminary notes}

Let $E$ and $F$ be a nonempty sets. We denote by $F^{E}$ the set of all the
mappings from $E$ to $F.$

Let $f:E\rightarrow F$ be a mapping and $B$ a subset of $F$ such that $%
f\left( E\right) \subset B.$ We denote by by $f\left\vert _{E}^{B}\right. $
the mapping $f\left\vert _{E}^{B}\right. :E\rightarrow B,$ $x\mapsto f(x).$

Let $\left( V,+,.\right) $ be a complex vector space and $W$ a nonempty
subset of $V.$ $W$ is said to be a convex cone if for every $x,y\in W$ and $%
\lambda \in 
\mathbb{R}
^{+\ast }$ we have $x+y,$ $\lambda x\in W$ .

For every Banach space $E$ we denote by $\mathcal{B}\left( E\right) $ the
Banach space of all bounded linear oprator from $E$ to itself.

For each $p\in ]1,+\infty \lbrack ,$we denote by $l_{p}\left( 
\mathbb{C}
\right) $ the set of all the sequences $z:=\left( z_{n}\right) _{n\in 
\mathbb{N}
^{\ast }},$ $z_{n}\in \mathbb{%
\mathbb{C}
}$ with $\underset{n=1}{\overset{+\infty }{\sum }}\left\vert
z_{n}\right\vert ^{P}<+\infty .$ $l_{p}\left( 
\mathbb{C}
\right) $ is a Banach space when it is endowed with the norm%
\begin{equation*}
\begin{array}{cccc}
\left\Vert .\right\Vert _{p}: & l_{p}\left( 
\mathbb{C}
\right) & \rightarrow & 
\mathbb{R}
\\ 
& z:=\left( z_{n}\right) _{n\in 
\mathbb{N}
^{\ast }} & \mapsto & \left\Vert z\right\Vert _{p}:=\left( \underset{n=1}{%
\overset{+\infty }{\sum }}\left\vert z_{n}\right\vert ^{p}\right) ^{\frac{1}{%
p}}%
\end{array}%
\end{equation*}%
we denote by $l_{\infty }\left( 
\mathbb{C}
\right) $ the set of all the bounded sequences $z:=\left( z_{n}\right)
_{n\in 
\mathbb{N}
^{\ast }},$ $z_{n}\in \mathbb{%
\mathbb{C}
}$. $l_{\infty }\left( 
\mathbb{C}
\right) $ is a Banach space when it is endowed with the norm%
\begin{equation*}
\begin{array}{cccc}
\left\Vert .\right\Vert _{\infty }: & l_{\infty }\left( 
\mathbb{C}
\right) & \rightarrow & 
\mathbb{R}
\\ 
& z:=\left( z_{n}\right) _{n\in 
\mathbb{N}
^{\ast }} & \mapsto & \left\Vert z\right\Vert _{\infty }:=\underset{n\in 
\mathbb{N}
^{\ast }}{\sup }\left( \left\vert z_{n}\right\vert \right)%
\end{array}%
\end{equation*}

Let $\mathfrak{A}$ be a C*-algebra. The identity element of $\mathfrak{A}$
with respect to the addition is then denoted by $0_{\mathfrak{A}}.$ The
C*-algebra $\mathfrak{A}$ is called unital if it has an identity element
with respect to the multiplication, that we will denote by $1_{\mathfrak{A}%
}, $ with the condition that $1_{\mathfrak{A}}\neq $ $0_{\mathfrak{A}}$.

An element $a\in \mathfrak{A}$\ is said to be self-adjoint if $a^{\ast }=a.$ 
$a$\ is said to be positive if $a^{\ast }=a$ and $\sigma (a)\subset 
\mathbb{R}
^{+}.$\textit{\ }

An element $x\in \mathfrak{A}$ is called strictly positive if for any
positive nonzero linear form $f$, we have $f\left( a\right) \in 
\mathbb{R}
^{+\ast }$.

\bigskip

\textbf{Proposition 2.1}

\textit{Let} $\mathfrak{A}$ \textit{be a unital C*-algebra and } $a\ $%
\textit{a self-adjoint element of }$\mathfrak{A}.$ \textit{Then }$\ a$ 
\textit{is a positive element of }$\mathfrak{A}$ \textit{if and only if}$\
a=b^{2}$ \textit{for some self-adjoint element }$b$\textit{\ of }$\mathfrak{A%
}.$

\bigskip

\textbf{Proposition 2.2.}

\textit{Let} $\mathfrak{A}$ \textit{be a unital C*-algebra} \textit{and }$%
a\in \mathfrak{A}$\textit{\ a positive element.} \textit{Then there exists a
unique element }$b\in \mathfrak{A}$ \textit{such that }$a=b^{2}.$ \textit{%
The element }$b$\textit{\ is then called the square root of }$a$\textit{\
and is denoted by }$a^{\frac{1}{2}}.$

\textit{\bigskip }

\textbf{Proposition 2.3.}

\textit{Let} $\mathfrak{A}$ \textit{be a unital C*-algebra} \textit{and }$%
a\in \mathfrak{A}$\textit{\ a positive element.} \textit{Then }$a$ \textit{%
is strictly positive if and only if }$a$ \textit{is invertible.}

\bigskip

The center $Z\left( \mathfrak{A}\right) $ of $\mathfrak{A}$ is the set of
all the elements $a\in \mathfrak{A}$ such that $ab=ba$ for every $b\in 
\mathfrak{A.}$ Each element of $Z\left( \mathfrak{A}\right) $ is called a
central element of $\mathfrak{A.}$

\bigskip

\textbf{Proposition 2.4.}

\textit{For every positif element }$a\in Z(\mathfrak{A)}$, $a^{\frac{1}{2}}$ 
\textit{belongs to} $Z(\mathfrak{A).}$

\bigskip

We consider the binary relation $\preccurlyeq $ defined on $\mathfrak{A}$ in
the following way%
\begin{equation*}
a\preccurlyeq b\ \text{if and only if }b-a\text{ is positive}
\end{equation*}%
Then $\left( \mathfrak{A,}\preccurlyeq \right) $\ is a partial order.

Let $\mathfrak{H}$ be a Hilbert $\mathfrak{A}$-module and $\left\langle
.,.\right\rangle $\ the inner product on $\mathfrak{H}.$ For every \ we set 
\begin{equation*}
\left\vert x\right\vert :=\left( \left\langle x,x\right\rangle \right) ^{%
\frac{1}{2}},\left\Vert x\right\Vert _{\mathfrak{H}}:=\left\Vert \left\vert
x\right\vert \right\Vert _{\mathfrak{A}},\text{ }x\in \mathfrak{H}
\end{equation*}

\bigskip

\textbf{Proposition 2.5. (\cite{LANC}, }page 6\textbf{)}

$\left\Vert .\right\Vert _{\mathfrak{H}}$\textit{\ is then a norm on }$%
\mathfrak{H}$\textit{\ which satisfies, for each\ }$x_{1},...,x_{N},$\textit{%
\ }$y_{1},...,y_{N}$\textit{\ }$\in $\textit{\ }$\mathfrak{H},$\textit{\ the
following property }%
\begin{equation*}
\left\Vert \underset{j=1}{\overset{N}{\sum }}\left\langle
x_{j},y_{j}\right\rangle \right\Vert _{\mathfrak{A}}^{2}\leq \left\Vert 
\underset{j=1}{\overset{N}{\sum }}\left\langle x_{j},x_{j}\right\rangle
\right\Vert _{\mathfrak{A}}\left\Vert \underset{j=1}{\overset{N}{\sum }}%
\left\langle y_{j},y_{j}\right\rangle \right\Vert _{\mathfrak{A}}
\end{equation*}

\bigskip

A Hilbert\textit{\ }$\mathfrak{A}$-module $\mathfrak{H}$ is said to be
countably generated if there exists a sequence $\left( x_{n}\right) _{n\in 
\mathbb{N}
^{\ast }}$ of elements of $\mathfrak{H}$ such that $\mathfrak{H}$ equals the
norm-closure of the $\mathfrak{A}$-linear span of the set $\left\{
x_{n}:n\in 
\mathbb{N}
^{\ast }\right\} .$

Let $\mathfrak{K}$ be a closed submodule of a Hilbert $\mathfrak{A}$-module $%
\mathfrak{H}$. The orthogonal complement $\mathfrak{K}^{\perp }$ of $%
\mathfrak{K}$ is the set 
\begin{equation*}
\mathfrak{K}^{\perp }:=\left\{ u\in \mathfrak{H:}\text{ }\left\langle
u,x\right\rangle =0_{\mathfrak{A}},\text{ }x\in \mathfrak{K}\right\}
\end{equation*}%
$\mathfrak{K}^{\perp }$ is then a closed submodule of the Hilbert $\mathfrak{%
A}$-module $\mathfrak{H.}$ We say that $\mathfrak{K}$ is orthogonally
complemented if 
\begin{equation*}
\mathfrak{H:=K\oplus K}^{\perp }
\end{equation*}%
The projection onto $\mathfrak{K}$ related to the direct sum is then called
the orthogonal projection on $\mathfrak{K}$ and is represented by the
notation $P_{\mathfrak{H}}.$

Let\textit{\ }$\mathfrak{H}$ and $\mathfrak{K}$ be a Hilbert\textit{\ }$%
\mathfrak{A}$-module. A mapping $T:\mathfrak{H\rightarrow K}$ is called an
adjointable operator if there exists a mapping\textit{\ }$S:\mathfrak{%
K\rightarrow H}$ such that%
\begin{equation*}
\left\langle T\left( u\right) ,v\right\rangle =\left\langle u,S\left(
v\right) \right\rangle ,\text{ }u\in \mathfrak{H},\text{ }v\in \mathfrak{K}
\end{equation*}%
We denote by $Hom_{\mathfrak{A}}^{\ast }\left( \mathfrak{H},\mathfrak{K}%
\right) $ the set of all $T:$ is called an adjointable operator from $%
\mathfrak{H\ }$to $\mathfrak{K.}$

\bigskip

\textbf{Proposition 2.6.}

\textit{Let }$\mathfrak{H}$ \textit{and} $\mathfrak{K}$ \textit{be a Hilbert 
}$\mathfrak{A}$-\textit{modules. }

a. \textit{If a mapping }$T:\mathfrak{H\rightarrow K}$ \textit{is
adjointable operator then there exists a unique mapping }$S:\mathfrak{%
K\rightarrow H}$ \textit{such that}%
\begin{equation*}
\left\langle T\left( u\right) ,v\right\rangle =\left\langle u,S\left(
v\right) \right\rangle ,\text{ }u\in \mathfrak{H},\text{ }v\in \mathfrak{K}
\end{equation*}%
\textit{The mapping }$S$\textit{\ is then called the adjoint of the mapping }%
$T$\textit{\ and is denoted by }$T^{\ast }.$

b. \textit{Every adjointable operator }$T:\mathfrak{H\rightarrow K}$\textit{%
\ together with its adjoint }$T^{\ast }:\mathfrak{K\rightarrow H}$ \textit{%
are a linear mapping of }$\mathfrak{A}$\textit{-modules and a bounded linear
operators of Banach spaces.}

c. $Hom_{\mathfrak{A}}^{\ast }\left( \mathfrak{H},\mathfrak{K}\right) $ 
\textit{is a Banach space with the usual operator norm.}

d. $Hom_{\mathfrak{A}}^{\ast }\left( \mathfrak{H},\mathfrak{H}\right) $ 
\textit{is a C*-algebra.}

\bigskip

\bigskip

\textbf{Theorem 2.7. }(\textbf{\cite{ARAM}, }page 472)

\textit{Let }$\mathfrak{H}$\textit{\ be a Hilbert }$\mathfrak{A}$\textit{%
-module, and }$T\in \mathcal{B}\left( \mathfrak{H}\right) $\textit{\ such
that }$T^{\ast }=T$\textit{. The following statements are equivalent :}

i. \ $T$ \textit{is invertible}

ii. \textit{There exists a real constants }$m,M>0$\textit{\ such that }%
\begin{equation*}
m\left\Vert x\right\Vert _{\mathfrak{H}}\leq \left\Vert T\left( x\right)
\right\Vert _{\mathfrak{H}}\leq M\left\Vert x\right\Vert _{\mathfrak{H}},%
\text{ }x\in \mathfrak{H}
\end{equation*}

ii. \textit{There exists a real constants }$m_{1},M_{1}>0$\textit{\ such that%
}%
\begin{equation*}
m_{1}\left\langle x,x\right\rangle _{\mathfrak{H}}\leq \left\langle T\left(
x\right) ,T(x)\right\rangle \leq M_{1}\left\langle x,x\right\rangle ,\text{ }%
x\in \mathfrak{H}
\end{equation*}

\bigskip

\bigskip

Let $\mathfrak{I}$\ an ideal of $\mathfrak{A}.$ $\mathfrak{I}$\ is said to
be an essential ideal of $A$ if the following implication holds for every
ideal $\mathfrak{N}$ of $\mathfrak{A}$%
\begin{equation*}
\mathfrak{N\cap I}=\{0_{\mathfrak{A}}\}\implies \mathfrak{N}=\{0_{\mathfrak{A%
}}\}
\end{equation*}

\bigskip

\bigskip

\textbf{Theorem 2.8. ((}\cite{LANC}), pages 14\textbf{)}

1. \textit{For any C*-algebra }$\mathfrak{A}$\textit{\ there exists a unique
(up to isomorphism) C*-algebra }$M\left( \mathfrak{A}\right) $\textit{\ such
that }$\mathfrak{A}$\textit{\ is an essential ideal of }$M\left( \mathfrak{A}%
\right) $\textit{\ and for every C*-algebra }$\mathcal{C}$\textit{\
containing }$\mathfrak{A}$\textit{\ is an essential ideal of }$\mathcal{C}$%
\textit{, there is a }$\mathfrak{A}$\textit{-homomorphism }$\Phi :\mathcal{C}
$\textit{\ }$\rightarrow $\textit{\ }$M\left( \mathfrak{A}\right) $\textit{\
of C*-algebras such that the restriction }$\Phi \left\vert _{\mathcal{C}%
}^{\Phi \left( \mathcal{C}\right) }\right. $\textit{\ is an }$\mathfrak{A}$%
\textit{-isomorphism of \ C*-algebras. }$M\left( \mathfrak{A}\right) $ 
\textit{is called the multiplier algebra of }$\mathfrak{A.}$

\textit{2. If }$\mathfrak{A}$\textit{\ is unital then }$M\left( \mathfrak{A}%
\right) =\mathfrak{A}.$

\bigskip

\bigskip

Let\textit{\ }$\mathfrak{A}$\textit{\ }be a C*-algebra and $\mathfrak{H}$ a
Hilbert $\mathfrak{A}$-module. We denote by\textit{\ }$\mathbb{J}_{\mathfrak{%
H}}$ the ideal of\textit{\ }$\mathfrak{A}$\ generated by\textit{\ }$%
\{\left\langle x,y\right\rangle _{E}:x,y\in \mathfrak{H}\}.$ We say that $%
\mathfrak{H}$\ is full if $\mathfrak{A}=\mathbb{J}_{\mathfrak{H}}.$

Let $\mathfrak{H}$\ and $\mathfrak{K}$\ be a Hilbert $\mathfrak{A}$-modules
and $\Psi :\mathfrak{H}$\ $\rightarrow $\ $\mathfrak{K}$\ a $\mathfrak{A}$%
-linear mapping (not necessarily bounded). We say that $\Psi $\ is
orthogonality preserving if 
\begin{equation*}
\left\langle x,y\right\rangle _{E}=0_{E}\implies \left\langle \Psi \left(
x\right) ,\Psi \left( y\right) \right\rangle _{F}=0_{F}\text{ }
\end{equation*}%
for every $x,y\in \mathfrak{H.}$

\bigskip

Combining the theorem 2.8.\ with the main result in (\cite{LEUN}), we obtain
the following result.

\bigskip

\textbf{Theorem 2.9. }

\textit{Let }$\mathfrak{H}$\textit{\ and }$\mathfrak{K}$\textit{\ be a
Hilbert} $\mathfrak{A}$\textit{-modules and }$\Psi :\mathfrak{H}$\textit{\ }$%
\rightarrow $\textit{\ }$\mathfrak{K}$\textit{\ a }$\mathfrak{A}$\textit{%
-linear mapping. We assume that }$\mathfrak{A}$\textit{\ is unital and that }%
$\mathfrak{H}$\textit{\ is full}.

a. $\Psi $\textit{\ is orthogonality preserving if and only if there exists
a central positive element }$\nu $ \textit{of} $\mathfrak{A}$ \textit{such
that}%
\begin{equation}
\left\langle \Psi \left( x\right) ,\Psi \left( y\right) \right\rangle
_{F}=\nu \left\langle x,y\right\rangle _{E},\text{ }x,y\in \mathfrak{H}
\label{ORTHOG}
\end{equation}

b. \textit{Assume that }$\Psi $\textit{\ is orthogonality preserving. Then }$%
\Psi $ \textit{is bounded and the central positive element }$\nu $ \textit{%
of }$\mathfrak{A}$ \textit{satisfying the relation (\ref{ORTHOG}) is unique.
So we denote }$\nu $ \textit{by} $\nu \left( \Psi \right) .$

c. \textit{Assume that }$\Psi $\textit{\ is bijective and orthogonality
preserving, then }$\nu \left( \Psi \right) ^{\frac{1}{2}}$\textit{\ is a
strictly positive element of} $\mathfrak{A}$\textit{\ invertible and }$\Psi :%
\mathfrak{H}$\textit{\ }$\rightarrow $\textit{\ }$\mathfrak{K}$ \textit{is
an isomorphism of Hilbert} $\mathfrak{A}$\textit{-module}.

\bigskip

\textbf{Proposition 2.10.}

a. \textit{Let }$\left( E,\left\Vert .\right\Vert _{E}\right) $\textit{, }$%
\left( F,\left\Vert .\right\Vert _{F}\right) $\textit{\ be Banach spaces. If 
}$\varphi :E\rightarrow F$\textit{\ is an isomorphism of Banach spaces, then}%
\begin{equation*}
\left\Vert \varphi ^{-1}\right\Vert _{\mathcal{B}\left( E,F\right)
}^{-1}\left\Vert x\right\Vert _{E}\leq \left\Vert \varphi \left( x\right)
\right\Vert _{F}\leq \left\Vert \varphi \right\Vert _{\mathcal{B}\left(
E,F\right) }\left\Vert x\right\Vert _{E},\text{ }x\in E
\end{equation*}

b. \textit{Let }$\mu $\textit{\ be a positive element of a C*-algebra }$%
\mathfrak{A.}$

i. \textit{Assume that }$\mu $\textit{\ is central in }$\mathfrak{A}.$%
\textit{\ Then the mapping }

\begin{equation*}
\begin{array}{cccc}
L_{\mu }: & \mathfrak{A} & \rightarrow & \mathfrak{A} \\ 
& a & \mapsto & \mu a%
\end{array}%
\end{equation*}%
\textit{is increasing.}

ii. \textit{Assume now that }$\mu $\textit{\ is strictly positive in }$%
\mathfrak{A}.$\textit{\ Then the following properties hold}%
\begin{equation*}
\left\{ 
\begin{array}{c}
\left\Vert \mu ^{-1}\right\Vert _{\mathfrak{A}}^{-1}\left\Vert a\right\Vert
_{\mathfrak{A}}\leq \left\Vert L_{\mu }\left( a\right) \right\Vert _{%
\mathfrak{H}},\text{ }a\in \mathfrak{A} \\ 
\left\Vert \mu ^{-1}\right\Vert _{\mathfrak{A}}^{-1}\left\Vert x\right\Vert
_{\mathfrak{H}}\leq \left\Vert \mu x\right\Vert _{\mathfrak{H}}\leq
\left\Vert \mu \right\Vert _{\mathfrak{A}}\left\Vert x\right\Vert _{%
\mathfrak{H}},\text{ }x\in \mathfrak{H}%
\end{array}%
\right.
\end{equation*}

\bigskip

\textbf{Proof}

1. Let $x\in E.$ It is clear that 
\begin{equation*}
\left\Vert \varphi \left( x\right) \right\Vert _{F}\leq \left\Vert \varphi
\right\Vert _{\mathcal{B}\left( E,F\right) }\left\Vert x\right\Vert _{E}
\end{equation*}%
On the other hand, we have we can write for every $x\in \mathfrak{A}$ 
\begin{eqnarray*}
\left\Vert x\right\Vert _{E} &=&\left\Vert \varphi ^{-1}\left( \varphi
\left( x\right) \right) \right\Vert _{E} \\
&\leq &\left\Vert \varphi ^{-1}\right\Vert _{\mathcal{B}\left( F,E\right)
}\left\Vert \varphi \left( x\right) \right\Vert _{F}\text{ }
\end{eqnarray*}%
It follows that%
\begin{equation*}
\left\Vert \varphi ^{-1}\right\Vert _{\mathcal{B}\left( F,E\right)
}^{-1}\left\Vert x\right\Vert _{E}\leq \left\Vert \varphi \left( x\right)
\right\Vert _{F}
\end{equation*}%
Consequently we obtain%
\begin{equation*}
\left\Vert \varphi ^{-1}\right\Vert _{\mathcal{B}\left( E,F\right)
}^{-1}\left\Vert x\right\Vert _{E}\leq \left\Vert \varphi \left( x\right)
\right\Vert _{F}\leq \left\Vert \varphi \right\Vert _{\mathcal{B}\left(
E,F\right) }\left\Vert x\right\Vert _{E},\text{ }x\in E
\end{equation*}

2. i. Assume that $\mu $ is a positive element of $\mathfrak{A.}$ Let $%
v,w\in \mathfrak{A}$ be such that $v\preccurlyeq w.$ Hence $w-v$ is a
positive element of $\mathfrak{A}$. We can then write%
\begin{eqnarray*}
L_{\mu }\left( w\right) -L_{\mu }\left( v\right) &=&w\mu -v\mu \\
&=&\left( w-v\right) \mu \\
&=&\left( w-v\right) ^{\frac{1}{2}}\left( w-v\right) ^{\frac{1}{2}}\mu \\
&=&\left( w-v\right) ^{\frac{1}{2}}\mu \left( w-v\right) ^{\frac{1}{2}} \\
&=&\left( \left( w-v\right) ^{\frac{1}{2}}\mu ^{\frac{1}{2}}\right) \left(
\mu ^{\frac{1}{2}}\left( w-v\right) ^{\frac{1}{2}}\right)
\end{eqnarray*}%
But the elements $\left( w-v\right) ^{\frac{1}{2}}$ and $\mu ^{\frac{1}{2}}$
are self-adjoint. It follows that%
\begin{equation*}
L_{\mu }\left( w\right) -L_{\mu }\left( v\right) =\left( \left( w-v\right) ^{%
\frac{1}{2}}\mu ^{\frac{1}{2}}\right) \left( \left( w-v\right) ^{\frac{1}{2}%
}\mu ^{\frac{1}{2}}\right) ^{\ast }
\end{equation*}%
Hence $L_{\mu }\left( w\right) -L_{\mu }\left( v\right) $ is a positive
element of $\mathfrak{A}$. It follows that $L_{\mu }\left( w\right)
\preccurlyeq L_{\mu }\left( v\right) $. Consequently the mapping $L_{\mu }$
is increasing.

ii. We have for every $a\in \mathfrak{A}$%
\begin{eqnarray*}
\left\Vert a\right\Vert _{\mathfrak{H}} &=&\left\Vert \mu ^{-1}\left( \mu
a\right) \right\Vert _{\mathfrak{H}} \\
&\leq &\left\Vert \mu ^{-1}\right\Vert _{\mathfrak{A}}\left\Vert L_{\mu
}\left( a\right) \right\Vert _{\mathfrak{A}}\text{ }
\end{eqnarray*}%
It follows that%
\begin{equation*}
\left\Vert \mu ^{-1}\right\Vert _{\mathfrak{A}}^{-1}\left\Vert a\right\Vert
_{\mathfrak{H}}\leq \left\Vert L_{\mu }\left( a\right) \right\Vert _{%
\mathfrak{A}},\text{ }a\in \mathfrak{A}
\end{equation*}

We have also for every $x\in \mathfrak{H}$%
\begin{eqnarray*}
\left\Vert \mu x\right\Vert _{\mathfrak{H}} &=&\sqrt{\left\Vert \left\langle
\mu x,\mu x\right\rangle \right\Vert _{\mathfrak{A}}} \\
&=&\sqrt{\left\Vert \mu \left\langle x,x\right\rangle \mu ^{\ast
}\right\Vert _{\mathfrak{A}}} \\
&\leq &\sqrt{\left\Vert \mu \right\Vert _{\mathfrak{A}}\left\Vert
\left\langle x,x\right\rangle \right\Vert _{\mathfrak{A}}\left\Vert \mu
^{\ast }\right\Vert _{\mathfrak{A}}} \\
&\leq &\left\Vert \mu \right\Vert _{\mathfrak{A}}\left\Vert x\right\Vert _{%
\mathfrak{H}}
\end{eqnarray*}%
Since $\mu $\ is strictly positive in $\mathfrak{A,}$ we can write for every 
$x\in \mathfrak{A}$ 
\begin{eqnarray*}
\left\Vert x\right\Vert _{\mathfrak{H}} &=&\left\Vert \mu ^{-1}\left( \mu
x\right) \right\Vert _{\mathfrak{H}} \\
&\leq &\left\Vert \mu ^{-1}\right\Vert _{\mathfrak{A}}\left\Vert \mu
x\right\Vert _{\mathfrak{H}}\text{ }
\end{eqnarray*}%
It follows that%
\begin{equation*}
\left\Vert \mu ^{-1}\right\Vert _{\mathfrak{A}}^{-1}\left\Vert x\right\Vert
_{\mathfrak{H}}\leq \left\Vert \mu x\right\Vert _{\mathfrak{H}}
\end{equation*}%
The final conclusion is that%
\begin{equation*}
\left\Vert \mu ^{-1}\right\Vert _{\mathfrak{A}}^{-1}\left\Vert x\right\Vert
_{\mathfrak{H}}\leq \left\Vert \mu x\right\Vert _{\mathfrak{H}}\leq
\left\Vert \mu \right\Vert _{\mathfrak{A}}\left\Vert x\right\Vert _{%
\mathfrak{H}},\text{ }x\in \mathfrak{H}
\end{equation*}%
Thence we achieve the proof of the proposition. $\square $

\bigskip

\bigskip

\section{*-fusion frames in Hilbert $\mathfrak{A}$-modules}

\bigskip

\subsection{\textbf{Definitions}}

\bigskip

\textbf{Definition 3.1.}

\textit{Let }$\mathfrak{A}$ \textit{be a unital\ C*-algebra and }$\mathfrak{H%
}$ \textit{a Hilbert} $\mathfrak{A}$-\textit{module and }$\left( \mathfrak{H}%
_{n}\right) _{n\in 
\mathbb{N}
^{\ast }}$ \textit{be} \textit{a sequence of orthogonally complemented
submodules of }$\mathfrak{H.}$ \textit{For each }$n\in 
\mathbb{N}
^{\ast }$\textit{\ we denote by }$P_{\mathfrak{H}_{n}}$ \textit{the
orthogonal projection of }$\mathfrak{H}$ \textit{onto} $\mathfrak{H}_{n}.$ 
\textit{\ }

a. \textit{A sequence }$\left( \omega _{n}\right) _{n\in 
\mathbb{N}
^{\ast }}$\textit{\ of central and strictly positive elements of }$\mathfrak{%
A}$ \textit{is called a weight of the C*-algebra} $\mathfrak{A.}$\textit{\
We denote by }$\mathcal{W}\left( \mathfrak{A}\right) $ \textit{the set of
all the weights of the C*-algebra }$\mathfrak{A.}$ \textit{The sequence} $%
\left( \left( \mathfrak{H}_{n},\omega _{n}\right) \right) _{n\in 
\mathbb{N}
^{\ast }}$ \textit{is} \textit{called a weighted sequence of orthogonally
complemented submodules of a }$\mathfrak{H.}$

b. \textit{Let} $\left( \left( \mathfrak{H}_{n},\omega _{n}\right) \right)
_{n\in 
\mathbb{N}
^{\ast }}$ \textit{be a weighted sequence of orthogonally complemented
submodules of a }$\mathfrak{H.}$ \textit{We say that} $\mathcal{(}\left( 
\mathfrak{H}_{n},\omega _{n}\right) )_{n\in 
\mathbb{N}
^{\ast }}$ \textit{is a *-fusion frame of }$\mathfrak{H}$ \textit{if for
each }$x\in \mathfrak{H,}$ \textit{the series }$\sum \omega
_{n}^{2}\left\vert P_{\mathfrak{H}_{n}}\left( x\right) \right\vert ^{2}$%
\textit{\ converges in the norm }$\left\Vert .\right\Vert _{\mathfrak{H}}$ 
\textit{and} \textit{there exist two strictly positive elements }$A$ \textit{%
and} $B$ \textit{of }$\mathfrak{A}$ \textit{such that\ } 
\begin{equation}
\left\vert Ax\right\vert ^{2}\preccurlyeq \underset{n=1}{\overset{+\infty }{%
\sum }}\omega _{n}^{2}\left\vert P_{\mathfrak{H}_{n}}\left( x\right)
\right\vert ^{2}\preccurlyeq \left\vert Bx\right\vert ^{2},\text{ }x\in 
\mathfrak{H}  \label{b}
\end{equation}%
\textit{The elements }$A$\textit{\ and }$B$\textit{\ are then respectively
called a lower and an upper bounds of }$\mathcal{G}.$ \textit{If }$A=$%
\textit{\ }$B,$ \textit{then the *-fusion frame }$\mathcal{G}$\ \textit{is
said to be tight (or more precisely }$A$\textit{-tight).} \textit{If }$A=$%
\textit{\ }$B=1_{\mathfrak{A}},$ \textit{then the *-fusion frame }$\mathcal{G%
}$ \ \textit{is said to be a Parseval-*-fusion frame of }$\mathfrak{H.}$

c. \textit{The set of all the sequences} $\left( \omega _{n}\right) _{n\in 
\mathbb{N}
^{\ast }}$ $\in $\textit{\ }$\mathcal{W}\left( \mathfrak{A}\right) $ \textit{%
such that} $\left( \left( \mathfrak{H}_{n},\omega _{n}\right) \right) _{n\in 
\mathbb{N}
^{\ast }}$ \textit{is a *-fusion frame of }$\mathfrak{H}$ \textit{is called
the set of multipliers of the sequence }$\left( \mathfrak{H}_{n}\right)
_{n\in 
\mathbb{N}
^{\ast }}$ \textit{and is denoted by }$\mathfrak{m}\left( \left( \mathfrak{H}%
_{n}\right) _{n\in 
\mathbb{N}
^{\ast }}\right) .$

\bigskip

\bigskip

\textbf{Definition 3.2.}

\textit{Let }$\mathcal{G}:=\left( \left( \mathfrak{H}_{n},\omega _{n}\right)
\right) _{n\in 
\mathbb{N}
^{\ast }}$ \textit{be a *-fusion frame of }$\mathfrak{H.}$

a. \textit{The mapping }%
\begin{equation*}
\begin{array}{cccc}
\mathcal{S}_{\mathcal{G}}: & \mathfrak{H} & \rightarrow & \mathfrak{H} \\ 
& x & \mapsto & \underset{n=1}{\overset{+\infty }{\sum }}\omega _{n}^{2}P_{%
\mathfrak{H}_{n}}\left( x\right)%
\end{array}%
\end{equation*}%
\textit{(if it is well-defined)} \textit{is then called the *-fusion frame
analysis operator associated to }$\mathcal{G}.$

b. \textit{The mapping }%
\begin{equation*}
\begin{array}{cccc}
\mathcal{T}_{\mathcal{G}}: & \mathfrak{H} & \rightarrow & l_{2}\left( 
\mathfrak{H}\right) \\ 
& x & \mapsto & \left( \omega _{n}P_{\mathfrak{H}_{n}}\left( x\right)
\right) _{n\in 
\mathbb{N}
^{\ast }}%
\end{array}%
\end{equation*}%
\textit{is called the *-fusion frame synthesis operator associated to }$%
\mathcal{G}.$

\bigskip

\bigskip

\bigskip

\subsection{Main results}

\bigskip

\textbf{Theorem 3.3.}

\textit{Let }$\mathcal{G}:=\left( \left( \mathfrak{H}_{n},\omega _{n}\right)
\right) _{n\in 
\mathbb{N}
^{\ast }}$ \textit{be a *-fusion frame of }$\mathfrak{H.}$\textit{Then} $%
\mathcal{T}_{\mathcal{G}}:\mathfrak{H\rightarrow }l_{2}\left( \mathfrak{H}%
\right) $ \textit{is a well-defined bounded linear adjointable operator and
the adjoint of }$\mathcal{T}_{\mathcal{G}}$ \textit{is the operator} 
\begin{equation*}
\begin{array}{cccc}
\mathcal{T}_{\mathcal{G}}^{\ast }: & l_{2}\left( \mathfrak{H}\right) & 
\mathfrak{\rightarrow } & \mathfrak{H} \\ 
& y:=\left( y_{n}\right) _{n\in 
\mathbb{N}
^{\ast }} & \mapsto & \mathcal{T}_{\mathcal{G}}^{\ast }\left( y\right) :=%
\underset{n=1}{\overset{+\infty }{\sum }}\omega _{n}P_{\mathfrak{H}%
_{n}}\left( y_{n}\right)%
\end{array}%
\end{equation*}

\bigskip

\textbf{Proof}

\textbf{1. Proof that the mapping }$\mathcal{T}_{\mathcal{G}}$ \textbf{is a
well-defined bounded linear operator.}

Let $x\in \mathfrak{H.}$ The series $\sum \left\vert \omega _{j}P_{\mathfrak{%
H}_{j}}\left( x\right) \right\vert ^{2}$ is then convergent. Hence the
sequence $\left( \omega _{n}P_{\mathfrak{H}_{n}}\left( x\right) \right)
_{n\in 
\mathbb{N}
^{\ast }}$ belongs to $l_{2}\left( \mathfrak{H}\right) $ and we have 
\begin{eqnarray*}
\left\Vert \left( \omega _{n}P_{\mathfrak{H}_{n}}\left( x\right) \right)
_{n\in 
\mathbb{N}
^{\ast }}\right\Vert _{l_{2}\left( \mathfrak{H}\right) } &=&\left\Vert 
\underset{n=1}{\overset{+\infty }{\sum }}\left\vert \omega _{n}P_{\mathfrak{H%
}_{n}}\left( x\right) \right\vert ^{2}\right\Vert _{\mathfrak{A}}^{\frac{1}{2%
}} \\
&\leq &\left\Vert \left\vert Bx\right\vert ^{2}\right\Vert _{\mathfrak{A}}^{%
\frac{1}{2}} \\
&\leq &\left\Vert Bx\right\Vert _{\mathfrak{H}} \\
&\leq &\left\Vert B\right\Vert _{\mathfrak{A}}\left\Vert x\right\Vert _{%
\mathfrak{H}}
\end{eqnarray*}%
Thus the mapping $\mathcal{T}_{\mathcal{G}}$ is a well-defined bounded
linear operator.

\textbf{2. Proof that the operator }$\mathcal{T}_{\mathcal{G}}$ \textbf{is
adjointable and determination of its adjoint}

Let $x\in \mathfrak{H,}$ $y:=\left( y_{j}\right) _{j\in 
\mathbb{N}
^{\ast }}\in l_{2}\left( \mathfrak{H}\right) ,$ $n,m\in 
\mathbb{N}
^{\ast }.$ To make the things simple we set 
\begin{equation*}
z_{n,m}:=\underset{j=n}{\overset{n+m}{\sum }}\omega _{j}P_{\mathfrak{H}%
_{j}}\left( y_{j}\right)
\end{equation*}%
We have, relying on proposition 2.5. 
\begin{eqnarray*}
\left\Vert z_{n,m}\right\Vert _{\mathfrak{H}}^{4} &=&\left\Vert \left\langle
z_{n,m},\underset{j=n}{\overset{n+m}{\sum }}\omega _{j}P_{\mathfrak{H}%
_{j}}\left( y_{j}\right) \right\rangle \right\Vert _{\mathfrak{A}}^{2} \\
&=&\left\Vert \underset{j=n}{\overset{n+m}{\sum }}\left\langle
z_{n,m},\omega _{j}P_{\mathfrak{H}_{j}}\left( y_{j}\right) \right\rangle
\right\Vert _{\mathfrak{A}}^{2} \\
&=&\left\Vert \underset{j=n}{\overset{n+m}{\sum }}\left\langle \omega _{j}P_{%
\mathfrak{H}_{j}}\left( z_{n,m}\right) ,y_{j}\right\rangle \right\Vert _{%
\mathfrak{A}}^{2} \\
&\leq &\left\Vert \underset{j=n}{\overset{n+m}{\sum }}\left\langle \omega
_{j}P_{\mathfrak{H}_{j}}\left( z_{n,m}\right) ,\omega _{j}P_{\mathfrak{H}%
_{j}}\left( z_{n,m}\right) \right\rangle \right\Vert _{\mathfrak{A}%
}\left\Vert \underset{j=n}{\overset{n+m}{\sum }}\left\langle
y_{j},y_{j}\right\rangle \right\Vert _{\mathfrak{A}} \\
&\leq &\left\Vert \left\vert Bz_{n,m}\right\vert ^{2}\right\Vert _{\mathfrak{%
A}}\left\Vert \underset{j=n}{\overset{n+m}{\sum }}\left\langle
y_{j},y_{j}\right\rangle \right\Vert _{\mathfrak{A}} \\
&\leq &\left\Vert Bz_{n,m}\right\Vert _{\mathfrak{H}}^{2}\left\Vert \underset%
{j=n}{\overset{n+m}{\sum }}\left\langle y_{j},y_{j}\right\rangle \right\Vert
_{\mathfrak{A}} \\
&\leq &\left\Vert B\right\Vert _{\mathfrak{A}}^{2}\left\Vert
z_{n,m}\right\Vert _{\mathfrak{H}}^{2}\left\Vert \underset{j=n}{\overset{n+m}%
{\sum }}\left\langle y_{j},y_{j}\right\rangle \right\Vert _{\mathfrak{A}}
\end{eqnarray*}%
Consequently%
\begin{equation}
\left\Vert z_{n,m}\right\Vert _{\mathfrak{H}}\leq \left\Vert B\right\Vert _{%
\mathfrak{A}}\sqrt{\left\Vert \underset{j=n}{\overset{n+m}{\sum }}%
\left\langle y_{j},y_{j}\right\rangle \right\Vert _{\mathfrak{A}}}  \label{a}
\end{equation}%
But since the series\ $\sum \left\langle y_{n},y_{n}\right\rangle $ is
convergent, we have 
\begin{equation*}
\underset{n\rightarrow +\infty }{\lim }\underset{m\in 
\mathbb{N}
^{\ast }}{\sup }\left\Vert \underset{j=n}{\overset{n+m}{\sum }}\left\langle
y_{j},y_{j}\right\rangle \right\Vert _{\mathfrak{A}}=0
\end{equation*}%
Hence the series $\sum \omega _{n}P_{\mathfrak{H}_{n}}\left( y_{n}\right) $
is convergent$.$ Thus the mapping%
\begin{equation*}
\begin{array}{cccc}
\mathcal{K}: & l_{2}\left( \mathfrak{H}\right) & \rightarrow & \mathfrak{H}
\\ 
& y:=\left( y_{n}\right) _{n\in 
\mathbb{N}
^{\ast }} & \mapsto & \underset{n=1}{\overset{+\infty }{\sum }}\omega _{n}P_{%
\mathfrak{H}_{n}}\left( y_{n}\right)%
\end{array}%
\end{equation*}%
is well-defined. It is clear that $\mathcal{K}$ is also linear. When taking $%
n=1$ and tending $m$\ to infinity, the relation (\ref{a}) becomes 
\begin{equation*}
\left\Vert \mathcal{K}\left( y\right) \right\Vert _{l_{2}\left( \mathfrak{H}%
\right) }\leq \left\Vert B\right\Vert _{\mathfrak{A}}\left\Vert y\right\Vert
_{l_{2}\left( \mathfrak{H}\right) },\text{ }y\in l_{2}\left( \mathfrak{H}%
\right)
\end{equation*}%
Thus the mapping $\mathcal{K}$ is a bounded linear operator.

Now we have 
\begin{eqnarray*}
\left\langle \mathcal{T}_{\mathcal{G}}\left( x\right) ,y\right\rangle
_{l_{2}\left( \mathfrak{H}\right) } &=&\underset{n=1}{\overset{+\infty }{%
\sum }}\left\langle \omega _{n}P_{\mathfrak{H}_{n}}\left( x\right)
,y_{n}\right\rangle _{\mathfrak{H}} \\
&=&\underset{n=1}{\overset{+\infty }{\sum }}\left\langle \omega _{n}P_{%
\mathfrak{H}_{n}}\left( x\right) ,P_{\mathfrak{H}_{n}}\left( y_{n}\right)
\right\rangle _{\mathfrak{H}} \\
&=&\underset{n=1}{\overset{+\infty }{\sum }}\left\langle x,\omega _{n}P_{%
\mathfrak{H}_{n}}\left( y_{n}\right) \right\rangle _{\mathfrak{H}} \\
&=&\left\langle x,\underset{n=1}{\overset{+\infty }{\sum }}\omega _{n}P_{%
\mathfrak{H}_{n}}\left( y_{n}\right) \right\rangle _{\mathfrak{H}} \\
&=&\left\langle x,\mathcal{K}\left( y\right) \right\rangle _{\mathfrak{H}}
\end{eqnarray*}%
It follows that the operator bounded linear adjointable operator and the
adjoint of $T_{\mathcal{G}}$ is the operator $\mathcal{K}$. That is $T_{%
\mathcal{G}}^{\ast }=\mathcal{K}.$

The proof of the proposition is achieved.

$\blacksquare $

\bigskip

\textbf{Theorem 3.4.}

\textit{Let }$\mathcal{G}:=\left( \left( \mathfrak{H}_{n},\omega _{n}\right)
\right) _{n\in 
\mathbb{N}
^{\ast }}$ \textit{be a *-fusion frame of }$\mathfrak{H.}$\textit{Then} $%
\mathcal{S}_{\mathcal{G}}$ \textit{is a well-defined bounded linear positive
and self-adjoint and invertible operator. Furthermore the following
reconstruction formula holds}%
\begin{equation*}
x=\underset{n=1}{\overset{+\infty }{\sum }}\omega _{n}^{2}P_{\mathfrak{H}%
_{n}}\left( \mathcal{S}_{\mathcal{G}}^{-1}\left( x\right) \right) ,\text{ }%
x\in \mathfrak{H}
\end{equation*}

\bigskip

\textbf{Proof}

\textbf{1. Proof that }$\mathcal{S}_{\mathcal{G}}$\textbf{\ is a
well-defined positive bounded linear and self-adjoint operator. }

For every $x\in \mathfrak{H}$ we have $\left( \omega _{n}\pi _{\mathfrak{H}%
_{n}}\left( x\right) \right) _{n\in 
\mathbb{N}
^{\ast }}\in l_{2}\left( \mathfrak{H}\right) $. Hence the series $\sum
\omega _{n}^{2}P_{\mathfrak{H}_{n}}\left( x\right) $ is convergent for every 
$x\in \mathfrak{H}$ and we have

\begin{eqnarray*}
\mathcal{T}_{\mathcal{G}}^{\ast }\mathcal{T}_{\mathcal{G}}\left( x\right) &=&%
\mathcal{T}_{\mathcal{G}}^{\ast }\left( \left( \omega _{n}P_{\mathfrak{H}%
_{n}}\left( x\right) \right) _{n\in 
\mathbb{N}
^{\ast }}\right) \\
&=&\text{ }\underset{n=1}{\overset{+\infty }{\sum }}\omega _{n}^{2}P_{%
\mathfrak{H}_{n}}\left( x\right)
\end{eqnarray*}%
It follows that the mapping $\mathcal{S}_{\mathcal{G}}$ is well-defined and $%
\mathcal{S}_{\mathcal{G}}=\mathcal{T}_{\mathcal{G}}^{\ast }\mathcal{T}_{%
\mathcal{G}}.$ Hence $\mathcal{S}_{\mathcal{G}}$ is a positive bounded
linear and self-adjoint operator.

\textbf{2.} \textbf{Proof that the operator }$\mathcal{S}_{\mathcal{G}}$%
\textbf{\ is invertible}

Since $\mathcal{S}_{\mathcal{G}}$ is a positive bounded linear and
self-adjoint operator, the operator $\mathcal{S}_{\mathcal{G}}$ has a square
root $\sqrt{\mathcal{S}_{\mathcal{G}}}$which also a positive bounded linear
and self-adjoint operator. Hence we have for every $x\in \mathfrak{H}$%
\begin{eqnarray*}
\left\Vert \sqrt{\mathcal{S}_{\mathcal{G}}}\left( x\right) \right\Vert _{%
\mathfrak{H}}^{2} &=&\left\Vert \left\langle \sqrt{\mathcal{S}_{\mathcal{G}}}%
\left( x\right) ,\sqrt{\mathcal{S}_{\mathcal{G}}}\left( x\right)
\right\rangle \right\Vert _{\mathfrak{A}} \\
&=&\left\Vert \left\langle \mathcal{S}_{\mathcal{G}}\left( x\right)
,x\right\rangle \right\Vert _{\mathfrak{A}} \\
&=&\left\Vert \left\langle \underset{n=1}{\overset{+\infty }{\sum }}\omega
_{n}^{2}P_{\mathfrak{H}_{n}}\left( x\right) ,x\right\rangle \right\Vert _{%
\mathfrak{A}} \\
&=&\left\Vert \left\langle \underset{n=1}{\overset{+\infty }{\sum }}\omega
_{n}^{2}\left( P_{\mathfrak{H}_{n}}P_{\mathfrak{H}_{n}}\right) \left(
x\right) ,x\right\rangle \right\Vert _{\mathfrak{A}} \\
&=&\left\Vert \left\langle \underset{n=1}{\overset{+\infty }{\sum }}\omega
_{n}^{2}\pi _{\mathfrak{H}_{n}}\left( x\right) ,\pi _{\mathfrak{H}%
_{n}}\left( x\right) \right\rangle \right\Vert _{\mathfrak{A}} \\
&=&\left\Vert \underset{n=1}{\overset{+\infty }{\sum }}\omega
_{n}^{2}\left\vert P_{\mathfrak{H}_{n}}\left( x\right) \right\vert
^{2}\right\Vert _{\mathfrak{A}}
\end{eqnarray*}%
It follows that%
\begin{eqnarray}
\left\Vert \sqrt{\mathcal{S}_{\mathcal{G}}}\left( x\right) \right\Vert _{%
\mathfrak{H}} &\geq &\left( \left\Vert \underset{n=1}{\overset{+\infty }{%
\sum }}\omega _{n}^{2}\left\vert P_{\mathfrak{H}_{n}}\left( x\right)
\right\vert ^{2}\right\Vert _{\mathfrak{A}}\right) ^{\frac{1}{2}}
\label{INVERTIBLE} \\
&\geq &\left( \left\Vert \left\vert Ax\right\vert ^{2}\right\Vert _{%
\mathfrak{A}}\right) ^{\frac{1}{2}}  \notag \\
&\geq &\left\Vert Ax\right\Vert _{\mathfrak{H}}  \notag
\end{eqnarray}%
But we know that $A$ is a strictly positive element of $\mathfrak{A.}$ It
follows from the proposition 2.10. that 
\begin{equation}
\left\Vert Ax\right\Vert _{\mathfrak{H}}\geq \left\Vert A^{-1}\right\Vert _{%
\mathfrak{A}}^{-1}\left\Vert x\right\Vert _{\mathfrak{H}}  \label{yes}
\end{equation}%
When combining (\ref{INVERTIBLE}) with (\ref{yes}) we obtain%
\begin{equation*}
\left\Vert \sqrt{\mathcal{S}_{\mathcal{G}}}\left( x\right) \right\Vert _{%
\mathfrak{H}}\geq \left\Vert A^{-1}\right\Vert _{\mathfrak{A}%
}^{-1}\left\Vert x\right\Vert _{\mathfrak{H}}
\end{equation*}%
Hence $\sqrt{\mathcal{S}_{\mathcal{G}}}$ is self adjoint and bounded below.
It follows theorem 2.7. that the operator $\sqrt{\mathcal{S}_{\mathcal{G}}}$
is invertible. Consequently the operator $\mathcal{S}_{\mathcal{G}}=\left( 
\sqrt{\mathcal{S}_{\mathcal{G}}}\right) ^{2}$ is also invertible.

\textbf{3. Proof of the reconstruction formula}

The operator $\mathcal{S}_{\mathcal{G}}$ is invertible. Hence we have for
every $x\in \mathfrak{H}$%
\begin{eqnarray*}
x &=&\mathcal{S}_{\mathcal{G}}\left( \mathcal{S}_{\mathcal{G}}^{-1}\left(
x\right) \right) \\
&=&\text{ }\underset{n=1}{\overset{+\infty }{\sum }}\omega _{n}^{2}P_{%
\mathfrak{H}_{n}}\left( \mathcal{S}_{\mathcal{G}}^{-1}\left( x\right) \right)
\end{eqnarray*}%
So the reconstruction formula holds true.

The proof of the proposition is achieved.

$\blacksquare $

\bigskip

The following result is a direct consequence of the reconstruction formula.

\bigskip

\textbf{Corollary 3.5.}

\textit{Let }$\mathfrak{A}$\textit{\ be a unital C*-algebra, }$\mathfrak{H}$%
\textit{\ a Hilbert} $\mathfrak{A}$\textit{-modules which has a *-fusion
frame sequence of orthogonally complemented submodule of} $\mathfrak{H}$%
\textit{. Then }$\mathfrak{H}$ \textit{is countably generated.}

\bigskip

\textbf{Theorem 3.6.}

\textit{Let }$\mathfrak{A}$\textit{\ be a unital C*-algebra, }$\mathfrak{H}$%
\textit{\ a Hilbert} $\mathfrak{A}$\textit{-modules, }$\left( \mathfrak{H}%
_{n}\right) _{n\in 
\mathbb{N}
^{\ast }}$ \textit{a sequence of orthogonally complemented submodule of} $%
\mathfrak{H}$\textit{\ and} $\left( \omega _{n}\right) _{n\in 
\mathbb{N}
^{\ast }}$\textit{\ a weight of the C*-algebra} $\mathfrak{A.}$

$\left( \left( \mathfrak{H}_{n},\omega _{n}\right) \right) _{n\in 
\mathbb{N}
^{\ast }}$ \textit{is a *-fusion frame of }$\mathfrak{H}$ \textit{if and
only if the series }$\sum \omega _{n}^{2}\left\vert P_{\mathfrak{H}%
_{n}}\left( x\right) \right\vert ^{2}$\textit{\ is convergent for every }$%
x\in \mathfrak{H}$\textit{\ and there exists a real constants }$c,$ $d$ $>0$%
\textit{\ such that}%
\begin{equation}
c\left\Vert x\right\Vert _{\mathfrak{H}}^{2}\leq \left\Vert \underset{n=1}{%
\overset{+\infty }{\sum }}\omega _{n}^{2}\left\vert P_{\mathfrak{H}%
_{n}}\left( x\right) \right\vert ^{2}\right\Vert _{\mathfrak{A}}\leq
d\left\Vert x\right\Vert _{\mathfrak{H}}^{2},\text{ }x\in \mathfrak{H}
\label{frame}
\end{equation}

\bigskip

\textbf{Proof }

1. Assume now that $\left( \left( \mathfrak{H}_{n},\omega _{n}\right)
\right) _{n\in 
\mathbb{N}
^{\ast }}$ is a *-fusion frame. Then the series $\sum \omega
_{n}^{2}\left\vert P_{\mathfrak{H}_{n}}\left( x\right) \right\vert ^{2}$\ is
convergent for every $x\in \mathfrak{H}$\textit{\ }and\textit{\ }there exist
strictly positive elements $A$ and $B$ of $\mathfrak{A}$ such that%
\begin{equation*}
\left\vert Ax\right\vert ^{2}\preccurlyeq \underset{n=1}{\overset{+\infty }{%
\sum }}\omega _{n}^{2}\left\vert P_{\mathfrak{H}_{n}}\left( x\right)
\right\vert ^{2}\preccurlyeq \left\vert Bx\right\vert ^{2},\text{ }x\in 
\mathfrak{H}
\end{equation*}%
It follows that%
\begin{equation*}
\left\Vert \left\vert Ax\right\vert ^{2}\right\Vert _{\mathfrak{A}}\leq
\left\Vert \underset{n=1}{\overset{+\infty }{\sum }}\omega
_{n}^{2}\left\vert P_{\mathfrak{H}_{n}}\left( x\right) \right\vert
^{2}\right\Vert _{\mathfrak{A}}\leq \left\Vert \left\vert Bx\right\vert
^{2}\right\Vert _{\mathfrak{A}},\text{ }x\in \mathfrak{H}
\end{equation*}%
It follows from proposition 2.10. that 
\begin{equation*}
\left\Vert A^{-1}\right\Vert _{\mathfrak{A}}^{-2}\left\Vert x\right\Vert _{%
\mathfrak{H}}^{2}\leq \left\Vert \underset{n=1}{\overset{+\infty }{\sum }}%
\omega _{n}^{2}\left\vert P_{\mathfrak{H}_{n}}\left( x\right) \right\vert
^{2}\right\Vert _{\mathfrak{A}}\leq \left\Vert B\right\Vert _{\mathfrak{A}%
}^{2}\left\Vert x\right\Vert _{\mathfrak{H}}^{2},\text{ }x\in \mathfrak{H}
\end{equation*}

2. Assume now that the series $\sum \omega _{n}^{2}\left\vert P_{\mathfrak{H}%
_{n}}\left( x\right) \right\vert ^{2}$\ is convergent for every $x\in 
\mathfrak{H}$\textit{\ }and that\textit{\ }there\textit{\ }exists a real
constants $c,$ $d$ $>0$ such that 
\begin{equation*}
c\left\Vert x\right\Vert _{\mathfrak{H}}^{2}\leq \left\Vert \underset{n=1}{%
\overset{+\infty }{\sum }}\omega _{n}^{2}\left\vert P_{\mathfrak{H}%
_{n}}\left( x\right) \right\vert ^{2}\right\Vert _{\mathfrak{A}}\leq
d\left\Vert x\right\Vert _{\mathfrak{H}}^{2},\text{ }x\in \mathfrak{H}
\end{equation*}%
Let $n,$ $m\in 
\mathbb{N}
^{\ast },$ $x\in \mathfrak{H.}$ We set

\begin{equation*}
u_{n,m}:=\underset{j=n}{\overset{n+m}{\sum }}\omega _{j}^{2}P_{\mathfrak{H}%
_{j}}\left( x\right)
\end{equation*}%
Then, relying on proposition 2.5., we obtain 
\begin{eqnarray*}
\left\Vert u_{n,m}\right\Vert _{\mathfrak{H}}^{4} &=&\left\Vert \left\langle
u_{n,m},\underset{j=n}{\overset{n+m}{\sum }}\omega _{j}P_{\mathfrak{H}%
_{j}}\left( \omega _{j}P_{\mathfrak{H}_{j}}\left( x\right) \right)
\right\rangle \right\Vert _{\mathfrak{A}}^{2} \\
&=&\left\Vert \underset{j=n}{\overset{n+m}{\sum }}\left\langle
u_{n,m},\omega _{j}P_{\mathfrak{H}_{j}}\left( \omega _{j}P_{\mathfrak{H}%
_{j}}\left( x\right) \right) \right\rangle \right\Vert _{\mathfrak{A}}^{2} \\
&=&\left\Vert \underset{j=n}{\overset{n+m}{\sum }}\left\langle \omega _{j}P_{%
\mathfrak{H}_{j}}\left( u_{n,m}\right) ,\omega _{j}P_{\mathfrak{H}%
_{j}}\left( x\right) \right\rangle \right\Vert _{\mathfrak{A}}^{2} \\
&\leq &\left\Vert \underset{j=n}{\overset{n+m}{\sum }}\left\langle \omega
_{j}P_{\mathfrak{H}_{j}}\left( u_{n,m}\right) ,\omega _{j}P_{\mathfrak{H}%
_{j}}\left( u_{n,m}\right) \right\rangle \right\Vert _{\mathfrak{A}%
}\left\Vert \underset{j=n}{\overset{n+m}{\sum }}\left\langle \omega _{j}P_{%
\mathfrak{H}_{j}}\left( x\right) ,\omega _{j}P_{\mathfrak{H}_{j}}\left(
x\right) \right\rangle \right\Vert _{\mathfrak{A}} \\
&\leq &d\left\Vert u_{n,m}\right\Vert _{\mathfrak{H}}^{2}\left\Vert \underset%
{j=n}{\overset{n+m}{\sum }}\left\vert \omega _{j}P_{\mathfrak{H}_{j}}\left(
x\right) \right\vert ^{2}\right\Vert _{\mathfrak{A}} \\
&\leq &d\left\Vert u_{n,m}\right\Vert _{\mathfrak{H}}^{2}\left\Vert \underset%
{j=n}{\overset{n+m}{\sum }}\left\vert \omega _{j}P_{\mathfrak{H}_{j}}\left(
x\right) \right\vert ^{2}\right\Vert _{\mathfrak{A}} \\
&\leq &d\left\Vert u_{n,m}\right\Vert _{\mathfrak{H}}^{2}\left\Vert \underset%
{j=n}{\overset{n+m}{\sum }}\left\vert \omega _{j}P_{\mathfrak{H}_{j}}\left(
x\right) \right\vert ^{2}\right\Vert _{\mathfrak{A}}
\end{eqnarray*}%
Consequently%
\begin{equation}
\left\Vert u_{n,m}\right\Vert _{\mathfrak{H}}\leq \sqrt{d}\sqrt{\left\Vert 
\underset{j=n}{\overset{n+m}{\sum }}\left\vert \omega _{j}P_{\mathfrak{H}%
_{j}}\left( x\right) \right\vert ^{2}\right\Vert _{\mathfrak{A}}}  \label{BB}
\end{equation}%
But since the series\ $\underset{j=n}{\overset{n+m}{\sum }}\left\vert \omega
_{j}P_{\mathfrak{H}_{j}}\left( x\right) \right\vert ^{2}$ is convergent, we
have 
\begin{equation*}
\underset{n\rightarrow +\infty }{\lim }\underset{m\in 
\mathbb{N}
^{\ast }}{\sup }\left\Vert \underset{j=n}{\overset{n+m}{\sum }}\left\vert
\omega _{j}P_{\mathfrak{H}_{j}}\left( x\right) \right\vert ^{2}\right\Vert _{%
\mathfrak{A}}=0
\end{equation*}%
Hence the series $\sum \omega _{n}^{2}P_{\mathfrak{H}_{n}}\left( x\right) $
is convergent$.$ Thus the mapping%
\begin{equation*}
\begin{array}{cccc}
\mathcal{L}: & \mathfrak{H} & \rightarrow & \mathfrak{H} \\ 
& x & \mapsto & \underset{n=1}{\overset{+\infty }{\sum }}\omega _{n}^{2}P_{%
\mathfrak{H}_{n}}\left( x\right)%
\end{array}%
\end{equation*}%
is well-defined. It is clear that $\mathcal{L}$ is also linear. When taking $%
n=1$ and tending $m$\ to infinity, the relation (\ref{BB}) becomes 
\begin{equation*}
\left\Vert \mathcal{L}\left( y\right) \right\Vert _{l_{2}\left( \mathfrak{H}%
\right) }\leq \sqrt{d}\sqrt{\left\Vert \underset{j=1}{\overset{+\infty }{%
\sum }}\left\vert \omega _{j}P_{\mathfrak{H}_{j}}\left( x\right) \right\vert
^{2}\right\Vert _{\mathfrak{A}}},\text{ }y\in l_{2}\left( \mathfrak{H}\right)
\end{equation*}%
Consequently we have%
\begin{equation*}
\left\Vert \mathcal{L}\left( y\right) \right\Vert _{l_{2}\left( \mathfrak{H}%
\right) }\leq d\left\Vert x\right\Vert _{\mathfrak{A}},\text{ }y\in
l_{2}\left( \mathfrak{H}\right)
\end{equation*}%
Thus the linear mapping $\mathcal{L}$ is bounded.

On the other hand we have for every $x\in \mathfrak{H}$ and $n\in 
\mathbb{N}
^{\ast }$

\begin{eqnarray}
\underset{n=1}{\overset{n}{\sum }}\omega _{j}^{2}\left\vert P_{\mathfrak{H}%
_{j}}\left( x\right) \right\vert ^{2} &=&\underset{n=1}{\overset{n}{\sum }}%
\omega _{j}^{2}\left\langle P_{\mathfrak{H}_{j}}\left( x\right) ,P_{%
\mathfrak{H}_{j}}\left( x\right) \right\rangle  \label{BBB} \\
&=&\left\langle \underset{j=1}{\overset{n}{\sum }}\omega _{j}^{2}P_{%
\mathfrak{H}_{j}}\left( x\right) ,x\right\rangle  \notag
\end{eqnarray}%
When taking $n=1$ and tending $m$\ to infinity, the relation (\ref{BBB})
becomes%
\begin{equation}
\underset{n=1}{\overset{+\infty }{\sum }}\omega _{j}^{2}\left\vert P_{%
\mathfrak{H}_{j}}\left( x\right) \right\vert ^{2}=\left\langle \mathcal{L}%
\left( x\right) ,x\right\rangle ,\text{ }x\in \mathfrak{H}  \label{C}
\end{equation}%
It follows that%
\begin{equation}
c\left\Vert x\right\Vert _{\mathfrak{H}}^{2}\leq \left\Vert \left\langle 
\mathcal{L}\left( x\right) ,x\right\rangle \right\Vert _{\mathfrak{A}}\leq
d\left\Vert x\right\Vert _{\mathfrak{H}}^{2},\text{ }x\in \mathfrak{H}
\label{BBBB}
\end{equation}

We can easily prove that the operator $\mathcal{L}$ is self-adjoint and
positive. It follows that $\mathcal{L}$ has a square root $\sqrt{\mathcal{L}}
$ which is also self-adjoint$.$ Thence the relation (\ref{BBBB}) becomes%
\begin{equation*}
c\left\Vert x\right\Vert _{\mathfrak{H}}^{2}\leq \left\Vert \left\langle 
\sqrt{\mathcal{L}}\left( x\right) ,\sqrt{\mathcal{L}}\left( x\right)
\right\rangle \right\Vert _{\mathfrak{A}}\leq d\left\Vert x\right\Vert _{%
\mathfrak{H}}^{2},\text{ }x\in \mathfrak{H}
\end{equation*}%
That is 
\begin{equation}
\sqrt{c}\left\Vert x\right\Vert _{\mathfrak{H}}\leq \left\Vert \sqrt{%
\mathcal{L}}\left( x\right) \right\Vert _{\mathfrak{H}}\leq \sqrt{d}%
\left\Vert x\right\Vert _{\mathfrak{H}},\text{ }x\in \mathfrak{H}
\label{BBBBB}
\end{equation}%
Consequently, by virtue of theorem 2.7., there exist a real constants $%
c_{1}, $ $d_{1}>0$ such that%
\begin{equation*}
c_{1}\left\vert x\right\vert ^{2}\leq \left\langle \sqrt{\mathcal{L}}\left(
x\right) ,\sqrt{\mathcal{L}}\left( x\right) \right\rangle \leq
d_{1}\left\vert x\right\vert ^{2},\text{ }x\in \mathfrak{H}
\end{equation*}%
That is%
\begin{equation*}
c_{1}\left\vert x\right\vert ^{2}\leq \left\langle \mathcal{L}\left(
x\right) ,x\right\rangle \leq d_{1}\left\vert x\right\vert ^{2},\text{ }x\in 
\mathfrak{H}
\end{equation*}%
It follows from (\ref{C}) that%
\begin{equation*}
c_{1}\left\vert x\right\vert ^{2}\leq \underset{n=1}{\overset{+\infty }{\sum 
}}\omega _{j}^{2}\left\vert P_{\mathfrak{H}_{j}}\left( x\right) \right\vert
^{2}\leq d_{1}\left\vert x\right\vert ^{2},\text{ }x\in \mathfrak{H}
\end{equation*}%
So%
\begin{equation}
\left\vert \left( \sqrt{c_{1}}1_{\mathfrak{A}}\right) x\right\vert ^{2}\leq 
\underset{n=1}{\overset{+\infty }{\sum }}\omega _{j}^{2}\left\vert P_{%
\mathfrak{H}_{j}}\left( x\right) \right\vert ^{2}\leq \left\vert \left( 
\sqrt{d_{1}}1_{\mathfrak{A}}\right) x\right\vert ^{2},\text{ }x\in \mathfrak{%
H}  \label{D}
\end{equation}%
Since $\sqrt{c_{1}}1_{\mathfrak{A}}$ and $\sqrt{d_{1}}1_{\mathfrak{A}}$ are
strictly positive elements of $\mathfrak{A,}$ it follows from (\ref{D}) that 
$\left( \left( \mathfrak{H}_{n},\omega _{n}\right) \right) _{n\in 
\mathbb{N}
^{\ast }}$ is a *-fusion-frame.

$\blacksquare $

\bigskip

\bigskip

\textbf{Theorem 3.7.}

\textit{Let }$\mathfrak{A}$\textit{\ be a unital and full C*-algebra, }$E$%
\textit{\ and }$F$\textit{\ be a Hilbert} $\mathfrak{A}$\textit{-modules and 
}$\Psi :E$\textit{\ }$\rightarrow $\textit{\ }$F$\textit{\ a bijective }$%
\mathfrak{A}$\textit{-linear mapping which is orthogonality preserving.} 
\textit{Let }$\mathcal{F}:=\left( \left( V_{n},\omega _{n}\right) \right)
_{n\in 
\mathbb{N}
^{\ast }}$ \textit{be a *-fusion frame of }$E\ $\textit{of bounds }$A,B$%
\textit{\ }$\in \mathfrak{A}.$

1. $\left( \Psi \left( V_{n}\right) \right) _{n\in 
\mathbb{N}
^{\ast }}$\textit{\ is then a sequence of orthogonally complemented
subspaces of }$F.$ \textit{Furthermore the following relation holds}%
\begin{equation}
\pi _{\Psi \left( V_{n}\right) }=\Psi P_{V_{n}}\Psi ^{-1}  \label{projection}
\end{equation}

2. $\left( \left( \Psi \left( V_{n}\right) ,\omega _{n}\right) \right)
_{n\in 
\mathbb{N}
^{\ast }}$ \textit{is a *-fusion frame of }$F.$

\bigskip

\textbf{Proof}

1. For every $n\in 
\mathbb{N}
^{\ast },$ $V_{n}$\textit{\ }is an orthogonally complemented subspace of $E.$
Hence we have%
\begin{equation*}
E=V_{n}\oplus V_{n}^{\perp }
\end{equation*}%
Thanks to theorem 2.9., it follows from the assumptions on $\Psi $ that it
is an isomorphism of Hilbert $\mathfrak{A}$-module from $E$ into $F.$ It
follows then that%
\begin{equation*}
F=\Psi \left( V_{n}\right) \oplus \Psi \left( V_{n}^{\perp }\right) ,\text{ }%
n\in 
\mathbb{N}
^{\ast }
\end{equation*}%
Since $\Psi $ is orthogonality preserving, it follows that%
\begin{equation*}
\Psi \left( V_{n}^{\perp }\right) \subset \left( \Psi \left( V_{n}\right)
\right) ^{\perp },\text{ }n\in 
\mathbb{N}
^{\ast }
\end{equation*}%
Hence 
\begin{equation}
F=\Psi \left( V_{n}\right) +\left( \Psi \left( V_{n}\right) \right) ^{\perp
},\text{ }n\in 
\mathbb{N}
^{\ast }  \label{direct1}
\end{equation}%
Let $n\in 
\mathbb{N}
^{\ast }$ and $a\in \Psi \left( V_{n}\right) \cap \left( \Psi \left(
V_{n}\right) \right) ^{\perp }.$ There exists then $b\in V_{n}$ such that $%
a=\Psi \left( b\right) \in \left( \Psi \left( V_{n}\right) \right) ^{\perp
}. $ Hence%
\begin{equation*}
\left\langle \Psi \left( b\right) ,\Psi \left( b\right) \right\rangle =0_{%
\mathfrak{A}}
\end{equation*}%
But 
\begin{equation*}
\left\langle \Psi \left( b\right) ,\Psi \left( b\right) \right\rangle =\nu
\left( \Psi \right) \left\langle b,b\right\rangle
\end{equation*}%
It follows that 
\begin{equation*}
\nu \left( \Psi \right) \left\langle b,b\right\rangle =0_{\mathfrak{A}}
\end{equation*}%
But $\nu \left( \Psi \right) $ is invertible. So $\left\langle
b,b\right\rangle =0.$ Hence $b=0_{E}.$ Consequently we have%
\begin{equation}
\Psi \left( V_{n}\right) \cap \left( \Psi \left( V_{n}\right) \right)
^{\perp }=\{0_{E}\},\text{ }n\in 
\mathbb{N}
^{\ast }  \label{direct2}
\end{equation}%
It follows from (\ref{direct1}) and (\ref{direct2}) that 
\begin{equation}
F=\Psi \left( V_{n}\right) \oplus \left( \Psi \left( V_{n}\right) \right)
^{\perp },\text{ }n\in 
\mathbb{N}
^{\ast }  \label{direct3}
\end{equation}%
Hence $\left( \Psi \left( V_{n}\right) \right) _{n\in 
\mathbb{N}
^{\ast }}$\ is a sequence of orthogonally complemented subspaces of $F.$

It follows from the relation (\ref{direct3}) that%
\begin{equation*}
\Psi \left( V_{n}\right) =\Psi \left( P_{V_{n}}\left( E\right) \right) ,%
\text{ }F=\Psi \left( P_{V_{n}}\left( E\right) \right) \oplus \left( \Psi
\left( P_{V_{n}}\left( E\right) \right) \right) ^{\perp },\text{ }n\in 
\mathbb{N}
^{\ast }
\end{equation*}%
Hence%
\begin{equation*}
\Psi \left( V_{n}\right) =\left( \Psi P_{V_{n}}\Psi ^{-1}\right) \left(
E\right) ,\text{ }F=\left( \Psi P_{V_{n}}\Psi ^{-1}\right) \left( E\right)
\oplus \left( \left( \Psi P_{V_{n}}\Psi ^{-1}\right) \left( E\right) \right)
^{\perp },\text{ }n\in 
\mathbb{N}
^{\ast }
\end{equation*}%
It follows that%
\begin{equation*}
\pi _{\Psi \left( V_{n}\right) }=\Psi P_{V_{n}}\Psi ^{-1},\text{ }n\in 
\mathbb{N}
^{\ast }
\end{equation*}

2. For every $y$ $\in F$ and $n\in 
\mathbb{N}
^{\ast },$ we have%
\begin{eqnarray*}
\omega _{n}^{2}\left\vert P_{\Psi \left( V_{n}\right) }\left( y\right)
\right\vert ^{2} &=&\omega _{n}^{2}\left\vert \Psi P_{V_{n}}\Psi ^{-1}\left(
y\right) \right\vert ^{2} \\
&=&\omega _{n}^{2}\left\vert \Psi \left( P_{V_{n}}\left( \Psi ^{-1}\left(
y\right) \right) \right) \right\vert ^{2} \\
&=&\omega _{n}^{2}\left\langle \Psi \left( P_{V_{n}}\left( \Psi ^{-1}\left(
y\right) \right) \right) ,\Psi \left( P_{V_{n}}\left( \Psi ^{-1}\left(
y\right) \right) \right) \right\rangle \\
&=&\nu \left( \Psi \right) \omega _{n}^{2}\left\langle P_{V_{n}}\left( \Psi
^{-1}\left( y\right) \right) ,P_{V_{n}}\left( \Psi ^{-1}\left( y\right)
\right) \right\rangle \\
&=&\nu \left( \Psi \right) \omega _{n}^{2}\left\vert P_{V_{n}}\left( \Psi
^{-1}\left( y\right) \right) \right\vert ^{2}
\end{eqnarray*}%
So we obtain%
\begin{equation}
\omega _{n}^{2}\left\vert P_{\Psi \left( V_{n}\right) }\left( y\right)
\right\vert ^{2}=\nu \left( \Psi \right) \omega _{n}^{2}\left\vert
P_{V_{n}}\left( \Psi ^{-1}\left( y\right) \right) \right\vert ^{2},\text{ }%
n\in 
\mathbb{N}
^{\ast },y\in F  \label{fusion}
\end{equation}%
Since $\mathcal{F}$ is a *-fusion frame it follows that the series $\sum
\omega _{n}^{2}\left\vert P_{V_{n}}\left( \Psi ^{-1}\left( y\right) \right)
\right\vert ^{2}$ is convergent in $\mathfrak{A}$ for every $y\in F.$ Hence
the series $\sum \omega _{n}^{2}\left\vert P_{\Psi \left( V_{n}\right)
}\left( y\right) \right\vert ^{2}$ is also convergent in $\mathfrak{A}$ for
every $y\in F.$ By virtue of the relatio (\ref{fusion}) we can then write
for every $y\in F$%
\begin{equation*}
\underset{n=1}{\overset{+\infty }{\sum }}\omega _{n}^{2}\left\vert P_{\Psi
\left( V_{n}\right) }\left( y\right) \right\vert ^{2}=\nu \left( \Psi
\right) \underset{n=1}{\overset{+\infty }{\sum }}\omega _{n}^{2}\left\vert
P_{V_{n}}\left( \Psi ^{-1}\left( y\right) \right) \right\vert ^{2}
\end{equation*}%
But we know that%
\begin{equation*}
\left\vert B\Psi ^{-1}\left( y\right) \right\vert ^{2}\preccurlyeq \underset{%
n=1}{\overset{+\infty }{\sum }}\omega _{n}^{2}\left\vert P_{V_{n}}\left(
\Psi ^{-1}\left( y\right) \right) \right\vert ^{2}\preccurlyeq \left\vert
B\Psi ^{-1}\left( y\right) \right\vert ^{2}
\end{equation*}%
It follows then that 
\begin{equation*}
\nu \left( \Psi \right) \left\vert A\Psi ^{-1}\left( y\right) \right\vert
^{2}\preccurlyeq \underset{n=1}{\overset{+\infty }{\sum }}\omega
_{n}^{2}\left\vert P_{\Psi \left( V_{n}\right) }\left( y\right) \right\vert
^{2}\preccurlyeq \nu \left( \Psi \right) \left\vert B\Psi ^{-1}\left(
y\right) \right\vert ^{2}
\end{equation*}%
Since $\nu \left( \Psi \right) ^{\frac{1}{2}}$ is a central and strictly
positive element of $\mathfrak{A,}$ it follows that 
\begin{equation*}
\left\vert \nu \left( \Psi \right) ^{\frac{1}{2}}A\Psi ^{-1}\left( y\right)
\right\vert _{\mathfrak{A}}^{2}\leq \left\Vert \underset{n=1}{\overset{%
+\infty }{\sum }}\omega _{n}^{2}\left\vert P_{\Psi \left( V_{n}\right)
}\left( y\right) \right\vert ^{2}\right\Vert _{\mathfrak{A}}\leq \left\vert
\nu \left( \Psi \right) ^{\frac{1}{2}}B\Psi ^{-1}\left( y\right) \right\vert
_{\mathfrak{A}}^{2}
\end{equation*}%
But $A$ and $B$\ are strictly positive elements of $\mathfrak{A}$ and $\nu
\left( \Psi \right) ^{\frac{1}{2}}$ is a central and strictly positive
element of $\mathfrak{A.}$ Hence $\nu \left( \Psi \right) ^{\frac{1}{2}}A$
and $\nu \left( \Psi \right) ^{\frac{1}{2}}B$ are also strictly positive
elements of $\mathfrak{A.}$ Consequently $\left( \left( \Psi \left(
V_{n}\right) ,\omega _{n}\right) \right) _{n\in 
\mathbb{N}
^{\ast }}$ is a *-fusion frame of $F.$

$\blacksquare $

\bigskip

\textbf{Theorem 3.8. }

\textit{For every }$\left( \mathfrak{\alpha }_{n}\right) _{n\in 
\mathbb{N}
^{\ast }}$, $\left( \mathfrak{\beta }_{n}\right) _{n\in 
\mathbb{N}
^{\ast }}\in \mathfrak{m}\left( \left( \mathfrak{H}_{n}\right) _{n\in 
\mathbb{N}
^{\ast }}\right) $\textit{\ } \textit{the sequence} $\left( \mathfrak{\alpha 
}_{n}+\mathfrak{\beta }_{n}\right) _{n\in 
\mathbb{N}
^{\ast }}$ \textit{belongs to} $\mathfrak{m}\left( \left( \mathfrak{H}%
_{n}\right) _{n\in 
\mathbb{N}
^{\ast }}\right) .$ \textit{It follows that} $\mathfrak{m}\left( \left( 
\mathfrak{H}_{n}\right) _{n\in 
\mathbb{N}
^{\ast }}\right) $ \textit{is a convex cone of the} $%
\mathbb{C}
$\textit{-vector space} $\mathfrak{A}^{%
\mathbb{N}
^{\ast }}.$

\bigskip

\textbf{Proof }

For each $n\in 
\mathbb{N}
^{\ast }$ the element $\mathfrak{\alpha }_{n}^{-1}\mathfrak{\beta }_{n}$ is
positif. It follows that $1_{\mathfrak{A}}+\mathfrak{\alpha }_{n}^{-1}%
\mathfrak{\beta }_{n}$ is also positif. Hence we have $\sigma \left( 
\mathfrak{\alpha }_{n}^{-1}\mathfrak{\beta }_{n}\right) \subset 
\mathbb{R}
^{+}$. It follows that $-1\notin \sigma \left( \mathfrak{\alpha }_{n}^{-1}%
\mathfrak{\beta }_{n}\right) .$ So $1_{\mathfrak{A}}+\mathfrak{\alpha }%
_{n}^{-1}\mathfrak{\beta }_{n}$ is invertible in $\mathfrak{A}$. Hence the
element $\mathfrak{\alpha }_{n}+\mathfrak{\beta }_{n}$ is also invertible in 
$\mathfrak{A}$. So $\mathfrak{\alpha }_{n}+\mathfrak{\beta }_{n}$\ is a
strictly positive element of $\mathfrak{A.}$ It is trivial that $\mathfrak{%
\alpha }_{n}+\mathfrak{\beta }_{n}$ is a central element of $\mathfrak{A.}$
Consequently $\left( \mathfrak{\alpha }_{n}+\mathfrak{\beta }_{n}\right)
_{n\in 
\mathbb{N}
^{\ast }}$ $\in \mathcal{W}\left( \mathfrak{H}\right) .$ On the other hand
there exist a strictly positive elements $A,$ $B,$ $C,$ $D$ of $\mathfrak{A}$
such that the following relations hold for every $x\in \mathfrak{H}$ 
\begin{equation*}
\left\{ 
\begin{array}{c}
\left\vert Ax\right\vert ^{2}\preccurlyeq \underset{n=1}{\overset{+\infty }{%
\sum }}\mathfrak{\alpha }_{n}^{2}\left\vert P_{\mathfrak{H}_{n}}\left(
x\right) \right\vert ^{2}\preccurlyeq \left\vert Bx\right\vert ^{2} \\ 
\left\vert Cx\right\vert ^{2}\preccurlyeq \underset{n=1}{\overset{+\infty }{%
\sum }}\mathfrak{\beta }_{n}^{2}\left\vert P_{\mathfrak{H}_{n}}\left(
x\right) \right\vert ^{2}\preccurlyeq \left\vert Dx\right\vert ^{2}%
\end{array}%
\right.
\end{equation*}%
But since all the $\mathfrak{\alpha }_{n}$, $\mathfrak{\beta }_{n}$ are
central strictly positifs elements of $\mathfrak{A,}$ it follows that%
\begin{equation*}
\mathfrak{\alpha }_{n}^{2}+\mathfrak{\beta }_{n}^{2}\preccurlyeq \left( 
\mathfrak{\alpha }_{n}+\mathfrak{\beta }_{n}\right) ^{2}\preccurlyeq 2%
\mathfrak{\alpha }_{n}^{2}+\mathfrak{\beta }_{n}^{2},\text{ }n\in 
\mathbb{N}
^{\ast }
\end{equation*}%
and that the two series $\sum \mathfrak{\alpha }_{n}^{2}\left\vert P_{%
\mathfrak{H}_{n}}\left( x\right) \right\vert ^{2}$ and $\sum \mathfrak{\beta 
}_{n}^{2}\left\vert P_{\mathfrak{H}_{n}}\left( x\right) \right\vert ^{2}$
are convergent in $\mathfrak{A.}$ Consequently, for every $x\in \mathfrak{H}$%
, the series $\sum \left( \mathfrak{\alpha }_{n}+\mathfrak{\beta }%
_{n}\right) ^{2}\left\vert P_{\mathfrak{H}_{n}}\left( x\right) \right\vert
^{2}$ is convergent in $\mathfrak{A}$ and we have%
\begin{equation*}
\left\vert Ax\right\vert ^{2}\preccurlyeq \underset{n=1}{\overset{+\infty }{%
\sum }}\left( \mathfrak{\alpha }_{n}+\mathfrak{\beta }_{n}\right)
^{2}\left\vert P_{\mathfrak{H}_{n}}\left( x\right) \right\vert ^{2}
\end{equation*}%
in addition to

\begin{equation*}
\begin{array}{c}
\underset{n=1}{\overset{+\infty }{\sum }}\left( \mathfrak{\alpha }_{n}+%
\mathfrak{\beta }_{n}\right) ^{2}\left\vert P_{\mathfrak{H}_{n}}\left(
x\right) \right\vert ^{2}\preccurlyeq \underset{n=1}{\overset{+\infty }{\sum 
}}2\left( \mathfrak{\alpha }_{n}^{2}+\mathfrak{\beta }_{n}^{2}\right)
\left\vert \pi _{\mathfrak{H}_{n}}\left( x\right) \right\vert ^{2} \\ 
\preccurlyeq 2\left( \left\vert Cx\right\vert ^{2}+\left\vert Dx\right\vert
^{2}\right) \\ 
\preccurlyeq 2C\left\langle x,x\right\rangle C^{\ast }+2D\left\langle
x,x\right\rangle D^{\ast } \\ 
\preccurlyeq 2C\left\langle x,x\right\rangle C+2D\left\langle
x,x\right\rangle D%
\end{array}%
\end{equation*}%
But we know, thanks to proposition 2.10. that%
\begin{eqnarray*}
\left\Vert A^{-1}\right\Vert _{\mathfrak{A}}^{-1}\left\Vert x\right\Vert _{%
\mathfrak{H}} &\leq &\left\Vert Ax\right\Vert _{\mathfrak{H}} \\
\left\Vert 2C\left\langle x,x\right\rangle C+2D\left\langle x,x\right\rangle
D\right\Vert _{\mathfrak{A}} &\leq &2\left( \left\Vert C\right\Vert _{%
\mathfrak{A}}^{2}+\left\Vert D\right\Vert _{\mathfrak{A}}^{2}\right)
\left\Vert x\right\Vert _{\mathfrak{H}}^{2}
\end{eqnarray*}%
It follows that%
\begin{equation*}
\left\{ 
\begin{array}{c}
\left\Vert A^{-1}\right\Vert _{\mathfrak{A}}^{-2}\left\Vert x\right\Vert _{%
\mathfrak{H}}^{2}\leq \left\Vert \underset{n=1}{\overset{+\infty }{\sum }}%
\left( \mathfrak{\alpha }_{n}+\mathfrak{\beta }_{n}\right) ^{2}\left\vert P_{%
\mathfrak{H}_{n}}\left( x\right) \right\vert ^{2}\right\Vert _{\mathfrak{A}}
\\ 
\left\Vert \underset{n=1}{\overset{+\infty }{\sum }}\left( \mathfrak{\alpha }%
_{n}+\mathfrak{\beta }_{n}\right) ^{2}\left\vert P_{\mathfrak{H}_{n}}\left(
x\right) \right\vert ^{2}\right\Vert _{\mathfrak{A}}\leq 2\left( \left\Vert
C\right\Vert _{\mathfrak{A}}^{2}+\left\Vert D\right\Vert _{\mathfrak{A}%
}^{2}\right) \left\Vert x\right\Vert _{\mathfrak{H}}^{2}%
\end{array}%
\right.
\end{equation*}%
Hencen, by virtue of theorem 3.6. that $\left( \left( \mathfrak{H}_{n},%
\mathfrak{\alpha }_{n}+\mathfrak{\beta }_{n}\right) \right) _{n\in 
\mathbb{N}
^{\ast }}$ is a *-fusion frame $\mathfrak{H.}$

Finally we conclude that $\left( \mathfrak{\alpha }_{n}+\mathfrak{\beta }%
_{n}\right) _{n\in 
\mathbb{N}
^{\ast }}\in \mathfrak{m}\left( \left( \mathfrak{H}_{n}\right) _{n\in 
\mathbb{N}
^{\ast }}\right) .$ It is also clear that $\left( \lambda \mathfrak{\alpha }%
_{n}\right) _{n\in 
\mathbb{N}
^{\ast }}\in \mathfrak{m}\left( \left( \mathfrak{H}_{n}\right) _{n\in 
\mathbb{N}
^{\ast }}\right) $ for every $\lambda \in 
\mathbb{R}
^{+\ast }.$ Consequently $\mathfrak{m}\left( \left( \mathfrak{H}_{n}\right)
_{n\in 
\mathbb{N}
^{\ast }}\right) $ is a convex cone of the $%
\mathbb{C}
$-vector space $\mathfrak{A}^{%
\mathbb{N}
^{\ast }}.$

$\blacksquare $

\bigskip

\section{\textbf{Examples of *-fusion frames }}

\bigskip

\textbf{Example 1}

We denote $\mathfrak{L}$ the $%
\mathbb{C}
$-vector space of all the sequences $z:=\left( z_{n}\right) _{n\in 
\mathbb{N}
^{\ast }},$ $z_{n}\in \mathbb{%
\mathbb{C}
}$ with 
\begin{equation*}
\underset{n=1}{\overset{+\infty }{\sum }}\frac{1}{2^{n}}\left\vert
z_{n}\right\vert ^{2}<+\infty
\end{equation*}

It is clear that $\mathfrak{L}$ is a commutative C*-algebra when it is
endowed with the mappings%
\begin{eqnarray*}
&&%
\begin{array}{cccc}
\cdot : & \mathfrak{L}^{2} & \rightarrow & \mathfrak{L} \\ 
& \left( \left( a_{n,1}\right) _{n\in 
\mathbb{N}
^{\ast }},\left( a_{n,2}\right) _{n\in 
\mathbb{N}
^{\ast }}\right) & \mapsto & \left( a_{n,1}\right) _{n\in 
\mathbb{N}
^{\ast }}\cdot \left( a_{n,2}\right) _{n\in 
\mathbb{N}
^{\ast }}:=\left( a_{n,1}a_{n,2}\right) _{n\in 
\mathbb{N}
^{\ast }}%
\end{array}
\\
&&%
\begin{array}{cccc}
^{\ast }: & \mathfrak{L} & \rightarrow & \mathfrak{L} \\ 
& \left( a_{n}\right) _{n\in 
\mathbb{N}
^{\ast }} & \mapsto & \left( \left( a_{n}\right) _{n\in 
\mathbb{N}
^{\ast }}\right) ^{\ast }:=\left( \overline{a_{n}}\right) _{n\in 
\mathbb{N}
^{\ast }}%
\end{array}%
\end{eqnarray*}%
Let us also observe that the C*-algebra $\mathfrak{L}$ is also unital, the
unit being the sequence $1_{\mathfrak{L}}:=\left( 1\right) _{n\in 
\mathbb{N}
^{\ast }}.$

On the other hand, $\mathfrak{L}$ is a $\mathfrak{L}$-module for the
well-defined operation%
\begin{equation*}
\begin{array}{cccc}
\circ : & \mathfrak{L}\times \mathfrak{L} & \rightarrow & \mathfrak{L} \\ 
& \left( \left( z_{n}\right) _{n\in 
\mathbb{N}
^{\ast }},\left( u_{n}\right) _{n\in 
\mathbb{N}
^{\ast }}\right) & \mapsto & \left( z_{n}\right) _{n\in 
\mathbb{N}
^{\ast }}\circ \left( u_{n}\right) _{n\in 
\mathbb{N}
^{\ast }}:=\left( z_{n}u_{n}\right) _{n\in 
\mathbb{N}
^{\ast }}%
\end{array}%
\end{equation*}%
We prove also easily that $\mathfrak{L}$ is a Hilbert $\mathfrak{L}$-module
if it is equipped with the well-defined $\mathfrak{L}$-inner product%
\begin{equation*}
\begin{array}{cccc}
\left\langle .,.\right\rangle : & \mathfrak{L}\times \mathfrak{L} & 
\rightarrow & \mathfrak{L} \\ 
& \left( \left( z_{n,1}\right) _{n\in 
\mathbb{N}
^{\ast }},\left( z_{n,2}\right) _{n\in 
\mathbb{N}
^{\ast }}\right) & \mapsto & \left\langle \left( z_{n,1}\right) _{n\in 
\mathbb{N}
^{\ast }},\left( z_{n,2}\right) _{n\in 
\mathbb{N}
^{\ast }}\right\rangle :=\left( z_{n}\overline{z_{n,2}}\right) _{n\in 
\mathbb{N}
^{\ast }}%
\end{array}%
\end{equation*}

Let $\left( I_{n}\right) _{n\in 
\mathbb{N}
^{\ast }}$ be a sequence of nonempty finite subsets of $%
\mathbb{N}
^{\ast }$ such that$\underset{n\in 
\mathbb{N}
^{\ast }}{\cup }I_{n}=%
\mathbb{N}
^{\ast }$. Let us consider, for every $k\in 
\mathbb{N}
^{\ast },$ the bounded linear operator%
\begin{equation*}
\begin{array}{cccc}
P_{n}: & \mathfrak{L} & \rightarrow & \mathfrak{L} \\ 
& \left( z_{n}\right) _{n\in 
\mathbb{N}
^{\ast }} & \mapsto & \left( d_{n,k}z_{n}\right) _{n\in 
\mathbb{N}
^{\ast }}%
\end{array}%
\end{equation*}%
where%
\begin{equation*}
d_{n,k}:=\left\{ 
\begin{array}{c}
1\text{ if }n\in I_{k} \\ 
0\text{ if }n\in 
\mathbb{N}
^{\ast }\backslash I_{k}%
\end{array}%
\right.
\end{equation*}

\bigskip

\textbf{Proposition 4.1.}

\textit{For every }$k\in 
\mathbb{N}
^{\ast },$\textit{\ the set }$U_{k}:=P_{k}\left( \mathfrak{L}\right) $%
\textit{\ is an orthogonally complemented submodule of }$\mathfrak{L.}$ 
\textit{Hence }$P_{k}$\textit{\ is, for every }$k\in 
\mathbb{N}
^{\ast },$\textit{\ the orthogonal projection of }$\mathfrak{L}$ \textit{onto%
} $U_{k}.$

\bigskip

\textbf{Proof}

Let $k\in 
\mathbb{N}
^{\ast }.$\ It is clear that $U_{k}$ is a closed submodule of $\mathfrak{L.}$
Let $\left( b_{n}\right) _{n\in 
\mathbb{N}
^{\ast }}\in U_{k}^{\perp }.$ Hence we have for every $j\in I_{k}$ 
\begin{eqnarray*}
&&\left( 0,0,...,0,\underset{j^{th}\text{ place}}{\underbrace{b_{j}},}%
0,0,...\right) \\
&=&\left\langle \left( b_{n}\right) _{n\in 
\mathbb{N}
^{\ast }},\left( 0,0,...,0,\underset{j^{th}\text{ place}}{\underbrace{v},}%
0,0,...\right) \right\rangle \\
&=&\left( 0,0,...,,0,...\right)
\end{eqnarray*}%
It follows that $U_{k}^{\perp }\subset Ker\left( P_{k}\right) .$ It is also
clear that $Ker\left( P_{k}\right) \subset U_{k}^{\perp }$ . Hence we have
for every $k\in \mathfrak{%
\mathbb{N}
}^{\ast },$ $U_{k}^{\perp }=Ker\left( P_{k}\right) $. We prove easily that $%
P_{k}\left( \mathfrak{L}\right) \oplus Ker\left( P_{k}\right) =\mathfrak{L.}$
Consequently 
\begin{equation*}
U_{k}\oplus U_{k}^{\perp }=\mathfrak{L}
\end{equation*}%
The final conclusion is that the set $U_{k}:=P_{k}\left( \mathfrak{L}\right) 
$, is for every $k\in 
\mathbb{N}
^{\ast },$\ an orthogonally complemented submodule of $\mathfrak{L}.$

$\blacksquare $

\bigskip

\textbf{Proposition 4.2.}

$\mathcal{W}\left( \mathfrak{L}\right) $ \textit{is the set of all the
sequences }$\left( \left( a_{n,k}\right) _{k\in 
\mathbb{N}
^{\ast }}\right) _{n\in 
\mathbb{N}
^{\ast }}$ of $\left( \mathfrak{L}\right) ^{%
\mathbb{N}
^{\ast }}$ \textit{such that}%
\begin{equation*}
\left\{ 
\begin{array}{c}
a_{n,k}\in 
\mathbb{R}
^{+\ast },\text{ }n,\text{ }k\in 
\mathbb{N}
^{\ast } \\ 
\overset{+\infty }{\underset{n=1}{\sum }}\frac{1}{2^{n}a_{n,k}^{2}}<+\infty ,%
\text{ }k\in 
\mathbb{N}
^{\ast }%
\end{array}%
\right.
\end{equation*}

\bigskip

\textbf{Proof}

A sequence $\left( \alpha _{n}\right) _{n\in 
\mathbb{N}
^{\ast }}\in \mathfrak{L}$ is a positive element of $\mathfrak{L}$ if and
only if there exists a sequence $\left( b_{n}\right) _{n\in 
\mathbb{N}
^{\ast }}\in \mathfrak{L}$ such that 
\begin{equation*}
\left( \alpha _{n}\right) _{n\in 
\mathbb{N}
^{\ast }}=\left( b_{n}\right) _{n\in 
\mathbb{N}
^{\ast }}\left( \left( b_{n}\right) _{n\in 
\mathbb{N}
^{\ast }}\right) ^{\ast }
\end{equation*}%
That is 
\begin{equation*}
\left\{ 
\begin{array}{c}
\alpha _{n}=\left\vert b_{n}\right\vert ^{2},\text{ }n\in 
\mathbb{N}
^{\ast } \\ 
\overset{+\infty }{\underset{n=1}{\sum }}\frac{\left\vert b_{n}\right\vert
^{2}}{2^{n}}<+\infty%
\end{array}%
\right.
\end{equation*}%
These conditions rewrite 
\begin{equation*}
\left\{ 
\begin{array}{c}
\alpha _{n}\in 
\mathbb{R}
^{+},\text{ }n\in 
\mathbb{N}
^{\ast } \\ 
\overset{+\infty }{\underset{n=1}{\sum }}\frac{\alpha _{n}}{2^{n}}<+\infty%
\end{array}%
\right.
\end{equation*}%
But since $\left( \alpha _{n}\right) _{n\in 
\mathbb{N}
^{\ast }}\in \mathfrak{L}$ it follows that $\alpha _{n}=\underset{%
n\rightarrow +\infty }{o}\left( 2^{\frac{n}{2}}\right) .$ Hence $\overset{%
+\infty }{\underset{n=1}{\sum }}\frac{\alpha _{n}}{2^{n}}<+\infty .$
Consequently $\left( \alpha _{n}\right) _{n\in 
\mathbb{N}
^{\ast }}\in \mathfrak{L}$ is a positive element of $\mathfrak{L}$ if and
only%
\begin{equation*}
\alpha _{n}\in 
\mathbb{R}
^{+},\text{ }n\in 
\mathbb{N}
^{\ast }
\end{equation*}

Now let $\left( \alpha _{n}\right) _{n\in 
\mathbb{N}
^{\ast }}\in \mathfrak{L}$ be a positive element of $\mathfrak{L.}$ $\left(
\alpha _{n}\right) _{n\in 
\mathbb{N}
^{\ast }}$ is invertible in $\mathfrak{L}$ if and only if there exists a
sequence $\left( c_{n}\right) _{n\in 
\mathbb{N}
^{\ast }}\in \mathfrak{L}$ such that 
\begin{equation*}
\left( \alpha _{n}\right) _{n\in 
\mathbb{N}
^{\ast }}\left( c_{n}\right) _{n\in 
\mathbb{N}
^{\ast }}=\left( 1\right) _{n\in 
\mathbb{N}
^{\ast }}
\end{equation*}%
That is 
\begin{equation*}
\left\{ 
\begin{array}{c}
\alpha _{n}c_{n}=1,\text{ }n\in 
\mathbb{N}
^{\ast } \\ 
\overset{+\infty }{\underset{n=1}{\sum }}\frac{1}{2^{n}}\left\vert
c_{n}\right\vert ^{2}<+\infty%
\end{array}%
\right.
\end{equation*}%
These conditions rewrite 
\begin{equation*}
\left\{ 
\begin{array}{c}
\alpha _{n}\in 
\mathbb{R}
^{+\ast },\text{ }n\in 
\mathbb{N}
^{\ast } \\ 
\overset{+\infty }{\underset{n=1}{\sum }}\frac{1}{2^{n}\alpha _{n}^{2}}%
<+\infty%
\end{array}%
\right.
\end{equation*}

The conclusion is that $\left( \alpha _{n}\right) _{n\in 
\mathbb{N}
^{\ast }}$ is a strictly positive element of $\mathfrak{L}$ if and only if%
\begin{equation*}
\left\{ 
\begin{array}{c}
\alpha _{n}\in 
\mathbb{R}
^{+\ast },\text{ }n\in 
\mathbb{N}
^{\ast } \\ 
\overset{+\infty }{\underset{n=1}{\sum }}\frac{1}{2^{n}\alpha _{n}^{2}}%
<+\infty%
\end{array}%
\right.
\end{equation*}

The final conclusion is that a sequence $\left( \left( a_{n,k}\right) _{k\in 
\mathbb{N}
^{\ast }}\right) _{n\in 
\mathbb{N}
^{\ast }}$ of $\left( \mathfrak{L}\right) ^{%
\mathbb{N}
^{\ast }}$ belongs to $\mathcal{W}\left( \mathfrak{L}\right) $ if and only
if 
\begin{equation*}
\left\{ 
\begin{array}{c}
a_{n,k}\in 
\mathbb{R}
^{+\ast },\text{ }n,\text{ }k\in 
\mathbb{N}
^{\ast } \\ 
\overset{+\infty }{\underset{n=1}{\sum }}\frac{1}{2^{n}a_{n,k}^{2}}<+\infty ,%
\text{ }k\in 
\mathbb{N}
^{\ast }%
\end{array}%
\right.
\end{equation*}

The proof of the proposition is then achieved.

$\blacksquare $

\bigskip

We denote by $\preccurlyeq ,$ as usual the partial order on the C*-algebra $%
\mathfrak{L}$ defined by positive elements of $\mathfrak{L.}$ In view of the
proof of the last proposition, it is then clear that for every $a:=\left(
a_{n}\right) _{n\in 
\mathbb{N}
^{\ast }},$ $b:=\left( b_{n}\right) _{n\in 
\mathbb{N}
^{\ast }}\in \mathfrak{L,}$ we have $a\preccurlyeq b$ if and only if $%
b_{n}-a_{n}\in 
\mathbb{R}
^{+}$ for every $n\in 
\mathbb{N}
^{\ast }$.

\bigskip

\textbf{Proposition 4.3.}

\textit{A sequence }$\left( a_{n}\right) _{n\in 
\mathbb{N}
^{\ast }}$ \textit{of elements of} $\mathcal{W}\left( \mathfrak{L}\right) $ 
\textit{belongs to} $\mathfrak{m}\left( \left( U_{n}\right) _{n\in 
\mathbb{N}
^{\ast }}\right) $ \textit{if and only if }%
\begin{equation*}
\left\{ 
\begin{array}{c}
0<\underset{n=1}{\overset{+\infty }{\sum }}a_{n,k}^{2}d_{n,k}<+\infty ,\text{
}k\in 
\mathbb{N}
^{\ast } \\ 
\underset{k\in 
\mathbb{N}
^{\ast }}{\sup }\left( 2^{k}\left( \underset{n=1}{\overset{+\infty }{\sum }}%
a_{n,k}^{2}d_{n,k}\right) \right) _{k\in 
\mathbb{N}
^{\ast }}<+\infty%
\end{array}%
\right.
\end{equation*}%
\textit{In this case }$\left( \left( U_{n}.,a_{n}\right) _{n\in 
\mathbb{N}
^{\ast }}\right) $\textit{\ is a }$\left( \sqrt{\underset{n=1}{\overset{%
+\infty }{\sum }}a_{n,k}^{2}d_{k,n}}\right) _{k\in 
\mathbb{N}
^{\ast }}$-\textit{tight *-fusion frame} $\left( \left( U_{n},a_{n}\right)
_{n\in 
\mathbb{N}
^{\ast }}\right) $\textit{\ of }$\mathfrak{L}.$

\bigskip

\textbf{Proof}

Let $\left( a_{n}\right) _{n\in 
\mathbb{N}
^{\ast }}$ $\in \mathcal{W}\left( \mathfrak{L}\right) .$

$\left( a_{n}\right) _{n\in 
\mathbb{N}
^{\ast }}:=\left( \left( a_{n,k}\right) _{k\in 
\mathbb{N}
^{\ast }}\right) _{n\in 
\mathbb{N}
^{\ast }}\in $ $\mathfrak{m}\left( \left( U_{n}\right) _{n\in 
\mathbb{N}
^{\ast }}\right) $ if and only if there exists two elements $\alpha :=\left(
\alpha _{n}\right) _{n\in 
\mathbb{N}
^{\ast }}$ and $\beta :=\left( \beta _{n}\right) _{n\in 
\mathbb{N}
^{\ast }}$ belonging to $\mathcal{W}\left( \mathfrak{L}\right) $ such that

i. for every $x:=\left( x_{n}\right) _{n\in 
\mathbb{N}
^{\ast }}$ $\in \mathfrak{L,}$ the series $\sum a_{n}^{2}\left\langle
P_{U_{n}}\left( x\right) ,P_{U_{n}}\left( x\right) \right\rangle $ is
convergent in $\mathfrak{L;}$

ii. the following relation holds%
\begin{equation*}
\alpha ^{2}\left\langle x,x\right\rangle \preccurlyeq \underset{n=1}{\overset%
{+\infty }{\sum }}a_{n}^{2}\left\langle P_{U_{n}}\left( x\right)
,P_{U_{n}}\left( x\right) \right\rangle \preccurlyeq \beta ^{2}\left\langle
x,x\right\rangle
\end{equation*}%
But%
\begin{equation*}
a_{n}^{2}\left\langle P_{U_{n}}\left( x\right) ,P_{U_{n}}\left( x\right)
\right\rangle =\left( a_{n,k}^{2}d_{n,k}\left\vert x_{k}\right\vert
^{2}\right) _{k\in 
\mathbb{N}
^{\ast }},\text{ }n\in 
\mathbb{N}
^{\ast }
\end{equation*}%
Hence $\left( a_{n}\right) _{n\in 
\mathbb{N}
^{\ast }}\in $ $\mathfrak{m}\left( \left( U_{n}\right) _{n\in 
\mathbb{N}
^{\ast }}\right) $ if and only if

i'. the series $\sum \left( a_{n,k}^{2}d_{n,k}\left\vert x_{k}\right\vert
^{2}\right) _{k\in 
\mathbb{N}
^{\ast }}$ is convergent in $\mathfrak{L}$ for every $x:=\left( x_{k}\right)
_{k\in 
\mathbb{N}
^{\ast }}\in \mathfrak{L;}$

ii'. $\left( \alpha _{k}^{2}\left\vert x_{k}\right\vert ^{2}\right) _{k\in 
\mathbb{N}
^{\ast }}\preccurlyeq \underset{n=1}{\overset{+\infty }{\sum }}\left(
a_{n,k}^{2}d_{k,n}\left\vert x_{k}\right\vert ^{2}\right) _{k\in 
\mathbb{N}
^{\ast }}\preccurlyeq \left( \beta _{k}^{2}\left\vert x_{k}\right\vert
^{2}\right) _{k\in 
\mathbb{N}
^{\ast }}.$

It follows then clearly that $\left( a_{n}\right) _{n\in 
\mathbb{N}
^{\ast }}\in $ $\mathfrak{m}\left( \left( U_{n}\right) _{n\in 
\mathbb{N}
^{\ast }}\right) $ if and only if

\begin{equation*}
\left\{ 
\begin{array}{c}
\underset{n=1}{\overset{+\infty }{\sum }}a_{n,k}^{2}d_{n,k}<+\infty ,k\in 
\mathbb{N}
^{\ast } \\ 
\underset{k=1}{\overset{+\infty }{\sum }}\frac{1}{2^{k}}\left( \left( 
\underset{n=1}{\overset{+\infty }{\sum }}a_{n,k}^{2}d_{n,k}\right)
\left\vert x_{k}\right\vert ^{2}\right) ^{2}<+\infty \\ 
\alpha _{k}^{2}\leq \underset{n=1}{\overset{+\infty }{\sum }}%
a_{n,k}^{2}d_{k,n}\leq \beta _{k}^{2}\text{, }k\in 
\mathbb{N}
^{\ast }%
\end{array}%
\right.
\end{equation*}%
holds for every $x:=\left( x_{n}\right) _{n\in 
\mathbb{N}
^{\ast }}$ $\in \mathfrak{L.}$ But since the mapping%
\begin{equation*}
\begin{array}{ccc}
\mathfrak{L} & \rightarrow & l_{2}\left( 
\mathbb{C}
\right) \\ 
\left( x_{n}\right) _{n\in 
\mathbb{N}
^{\ast }} & \mapsto & \left( \frac{1}{2^{n}}x_{n}^{2}\right) _{n\in 
\mathbb{N}
^{\ast }}%
\end{array}%
\end{equation*}%
is onto, we obtain the following equivalence :

$\left( a_{n}\right) _{n\in 
\mathbb{N}
^{\ast }}\in $ $\mathfrak{m}\left( \left( U_{n}\right) _{n\in 
\mathbb{N}
^{\ast }}\right) $ if and only if

\begin{equation*}
\left\{ 
\begin{array}{c}
\underset{n=1}{\overset{+\infty }{\sum }}a_{n,k}^{2}d_{n,k}<+\infty ,\text{ }%
k\in 
\mathbb{N}
^{\ast } \\ 
\underset{k=1}{\overset{+\infty }{\sum }}2^{k}\left( \left( \underset{n=1}{%
\overset{+\infty }{\sum }}a_{n,k}^{2}d_{n,k}\right) \right) ^{2}\left\vert
y_{k}\right\vert ^{2}<+\infty \\ 
\alpha _{k}^{2}\leq \underset{n=1}{\overset{+\infty }{\sum }}%
a_{n,k}^{2}d_{k,n}\leq \beta _{k}^{2}\text{, }k\in 
\mathbb{N}
^{\ast }%
\end{array}%
\right.
\end{equation*}%
holds for every $y:=\left( y_{n}\right) _{n\in 
\mathbb{N}
^{\ast }}$ $\in l_{1}\left( 
\mathbb{C}
\right) .$

Finally it follows easily that $\left( a_{n}\right) _{n\in 
\mathbb{N}
^{\ast }}\in $ $\mathfrak{m}\left( \left( U_{n}\right) _{n\in 
\mathbb{N}
^{\ast }}\right) $ if and only if 
\begin{equation*}
\left\{ 
\begin{array}{c}
\underset{n=1}{\overset{+\infty }{\sum }}a_{n,k}^{2}d_{n,k}<+\infty ,\text{ }%
k\in 
\mathbb{N}
^{\ast } \\ 
\underset{k\in 
\mathbb{N}
^{\ast }}{\sup }\left( 2^{k}\left( \underset{n=1}{\overset{+\infty }{\sum }}%
a_{n,k}^{2}d_{n,k}\right) \right) _{k\in 
\mathbb{N}
^{\ast }}<+\infty \\ 
0<\underset{n=1}{\overset{+\infty }{\sum }}a_{n,k}^{2}d_{k,n},\text{ }k\in 
\mathbb{N}
^{\ast } \\ 
\underset{k=1}{\overset{+\infty }{\sum }}\underset{n=1}{\overset{+\infty }{%
\sum }}\frac{a_{n,k}^{2}d_{k,n}}{2^{k}}<+\infty%
\end{array}%
\right.
\end{equation*}%
That is 
\begin{equation*}
\left\{ 
\begin{array}{c}
0<\underset{n=1}{\overset{+\infty }{\sum }}a_{n,k}^{2}d_{n,k}<+\infty ,\text{
}k\in 
\mathbb{N}
^{\ast } \\ 
\underset{k\in 
\mathbb{N}
^{\ast }}{\sup }\left( 2^{k}\left( \underset{n=1}{\overset{+\infty }{\sum }}%
a_{n,k}^{2}d_{n,k}\right) \right) _{k\in 
\mathbb{N}
^{\ast }}<+\infty%
\end{array}%
\right.
\end{equation*}%
In this case we can take for $\alpha :=\left( \alpha _{n}\right) _{n\in 
\mathbb{N}
^{\ast }}$ and $\beta :=\left( \beta _{n}\right) _{n\in 
\mathbb{N}
^{\ast }}$ the values

\begin{equation*}
\alpha =\beta =\left( \sqrt{\underset{n=1}{\overset{+\infty }{\sum }}%
a_{n,k}^{2}d_{k,n}}\right) _{k\in 
\mathbb{N}
^{\ast }}
\end{equation*}%
and $\left( \left( U_{n},a_{n}\right) \right) _{n\in 
\mathbb{N}
^{\ast }}$ becomes a $\left( \sqrt{\underset{n=1}{\overset{+\infty }{\sum }}%
a_{n,k}^{2}d_{k,n}}\right) _{k\in 
\mathbb{N}
^{\ast }}$-tight *-fusion frame of $\mathfrak{L.}$

The proof of the proposition is then complete.

$\blacksquare $

\bigskip

\textbf{Example 2}

Let $\mathbb{H}$ be the well-known field of quaternions. We denote by $%
l_{\infty }\left( \mathbb{H}\right) $ the $%
\mathbb{C}
$-vector space of all the sequences $z:=\left( z_{n}\right) _{n\in 
\mathbb{N}
^{\ast }},$ $z_{n}\in \mathbb{H}$ with 
\begin{equation*}
\underset{n\in 
\mathbb{N}
^{\ast }}{\sup }\left( \left\vert z_{n}\right\vert \right) <+\infty
\end{equation*}%
and by $\left( \mathbb{H}\right) _{l_{2}}$ the $%
\mathbb{C}
$-vector space of all the sequences $z:=\left( z_{n}\right) _{n\in 
\mathbb{N}
^{\ast }},$ $z_{n}\in \mathbb{H}$ with 
\begin{equation*}
\underset{n=1}{\overset{+\infty }{\sum }}\left\vert z_{n}\right\vert
^{2}<+\infty
\end{equation*}

It is clear that $l_{\infty }\left( \mathbb{H}\right) $ is a noncommutative
C*-algebra when it is endowed with the mappings%
\begin{eqnarray*}
&&%
\begin{array}{cccc}
\cdot : & l_{\infty }\left( \mathbb{H}\right) ^{2} & \rightarrow & l_{\infty
}\left( \mathbb{H}\right) \\ 
& \left( \left( a_{n,1}\right) _{n\in 
\mathbb{N}
^{\ast }},\left( a_{n,2}\right) _{n\in 
\mathbb{N}
^{\ast }}\right) & \mapsto & \left( a_{n,1}\right) _{n\in 
\mathbb{N}
^{\ast }}\cdot \left( a_{n,2}\right) _{n\in 
\mathbb{N}
^{\ast }}:=\left( a_{n,1}a_{n,2}\right) _{n\in 
\mathbb{N}
^{\ast }}%
\end{array}
\\
&&%
\begin{array}{cccc}
^{\ast }: & l_{\infty }\left( \mathbb{H}\right) & \rightarrow & l_{\infty
}\left( \mathbb{H}\right) \\ 
& \left( a_{n}\right) _{n\in 
\mathbb{N}
^{\ast }} & \mapsto & \left( \left( a_{n}\right) _{n\in 
\mathbb{N}
^{\ast }}\right) ^{\ast }:=\left( \overline{a_{n}}\right) _{n\in 
\mathbb{N}
^{\ast }}%
\end{array}%
\end{eqnarray*}%
Let us also observe that the C*-algebra $l_{\infty }\left( \mathbb{H}\right) 
$ is also unital, the unit being the sequence $1_{l_{\infty }\left( \mathbb{H%
}\right) }:=\left( 1\right) _{n\in 
\mathbb{N}
^{\ast }}.$ On the other hand, $\left( \mathbb{H}\right) _{l_{2}}$ is a
(left) $l_{\infty }\left( \mathbb{H}\right) $-module for the well-defined
operation%
\begin{equation*}
\begin{array}{cccc}
\circ : & l_{\infty }\left( \mathbb{H}\right) \times \left( \mathbb{H}%
\right) _{l_{2}} & \rightarrow & \left( \mathbb{H}\right) _{l_{2}} \\ 
& \left( \left( z_{n}\right) _{n\in 
\mathbb{N}
^{\ast }},\left( u_{n}\right) _{n\in 
\mathbb{N}
^{\ast }}\right) & \mapsto & \left( z_{n}\right) _{n\in 
\mathbb{N}
^{\ast }}\circ \left( u_{n}\right) _{n\in 
\mathbb{N}
^{\ast }}:=\left( z_{n}u_{n}\right) _{n\in 
\mathbb{N}
^{\ast }}%
\end{array}%
\end{equation*}%
We prove also easily that $\left( \mathbb{H}\right) _{l_{2}}$ is a (left)
Hilbert $l_{\infty }\left( \mathbb{H}\right) $-module if it is equipped with
the well-defined $l_{\infty }\left( \mathbb{H}\right) $-inner product%
\begin{equation*}
\begin{array}{cccc}
\left\langle .,.\right\rangle : & \left( \mathbb{H}\right) _{l_{2}}\times
\left( \mathbb{H}\right) _{l_{2}} & \rightarrow & l_{\infty }\left( \mathbb{H%
}\right) \\ 
& \left( \left( z_{n,1}\right) _{n\in 
\mathbb{N}
^{\ast }},\left( z_{n,2}\right) _{n\in 
\mathbb{N}
^{\ast }}\right) & \mapsto & \left\langle \left( z_{n,1}\right) _{n\in 
\mathbb{N}
^{\ast }},\left( z_{n,2}\right) _{n\in 
\mathbb{N}
^{\ast }}\right\rangle :=\left( z_{n}\overline{z_{n,2}}\right) _{n\in 
\mathbb{N}
^{\ast }}%
\end{array}%
\end{equation*}%
(The fact that the complex algebra $\left( \mathbb{H}\right) _{l_{2}}$ is
complete is based on a result stated in ((\cite{KOTE}); page 359).

Let $\left( I_{n}\right) _{n\in 
\mathbb{N}
^{\ast }}$ be a sequence of nonempty finite subsets of $%
\mathbb{N}
^{\ast }$ such that$\underset{n\in 
\mathbb{N}
^{\ast }}{\cup }I_{n}=%
\mathbb{N}
^{\ast }$. Let us consider, for every $k\in 
\mathbb{N}
^{\ast },$ the bouded linear operator%
\begin{equation*}
\begin{array}{cccc}
P_{k}: & \left( \mathbb{H}\right) _{l_{2}} & \rightarrow & \left( \mathbb{H}%
\right) _{l_{2}} \\ 
& \left( z_{n}\right) _{n\in 
\mathbb{N}
^{\ast }} & \mapsto & \left( d_{k,n}z_{n}\right) _{n\in 
\mathbb{N}
^{\ast }}%
\end{array}%
\end{equation*}%
where%
\begin{equation*}
d_{k,n}:=\left\{ 
\begin{array}{c}
1\text{ if }n\in I_{k} \\ 
0\text{ if }n\in 
\mathbb{N}
^{\ast }\backslash I_{k}%
\end{array}%
\right.
\end{equation*}

We obtain the following result in a similar as for the previous example.

\bigskip

\textbf{Proposition 4.4.}

\textit{For every }$k\in 
\mathbb{N}
^{\ast },$\textit{\ the set }$U_{k}:=P_{k}\left( \left( \mathbb{H}\right)
_{l_{2}}\right) $\textit{\ is an orthogonally complemented submodule of }$%
\left( \mathbb{H}\right) _{l_{2}}\mathfrak{.}$\textit{Hence }$P_{k}$\textit{%
\ is, for every }$k\in 
\mathbb{N}
^{\ast },$\textit{\ the orthogonal projection of }$\left( \mathbb{H}\right)
_{l_{2}}$ \textit{onto} $U_{k}.$

\bigskip

\textbf{Proposition 4.5.}

$\mathcal{W}\left( l_{\infty }\left( \mathbb{H}\right) \right) $ \textit{is
the set of all the sequences }$\left( \left( a_{n,k}\right) _{k\in 
\mathbb{N}
^{\ast }}\right) _{n\in 
\mathbb{N}
^{\ast }}$ of $\left( l_{\infty }\left( \mathbb{H}\right) \right) ^{%
\mathbb{N}
^{\ast }}$ \textit{such that}%
\begin{equation*}
\left\{ 
\begin{array}{c}
a_{n,k}\in 
\mathbb{R}
^{+\ast },\text{ }n,\text{ }k\in 
\mathbb{N}
^{\ast } \\ 
0<\underset{n\in 
\mathbb{N}
^{\ast }}{\inf }\left( a_{n,k}\right) ,\text{ }k\in 
\mathbb{N}
^{\ast }%
\end{array}%
\right.
\end{equation*}

\bigskip

\textbf{Proof}

A sequence $\left( \alpha _{n}\right) _{n\in 
\mathbb{N}
^{\ast }}\in l_{\infty }\left( \mathbb{H}\right) $ is a positive element of $%
l_{\infty }\left( \mathbb{H}\right) $ if and only if there exists a sequence 
$\left( b_{n}\right) _{n\in 
\mathbb{N}
^{\ast }}\in l_{\infty }\left( \mathbb{H}\right) $ such that 
\begin{equation*}
\left( \alpha _{n}\right) _{n\in 
\mathbb{N}
^{\ast }}=\left( b_{n}\right) _{n\in 
\mathbb{N}
^{\ast }}\left( \left( b_{n}\right) _{n\in 
\mathbb{N}
^{\ast }}\right) ^{\ast }
\end{equation*}%
That is 
\begin{equation*}
\left\{ 
\begin{array}{c}
\alpha _{n}=\left\vert b_{n}\right\vert ^{2},\text{ }n\in 
\mathbb{N}
^{\ast } \\ 
\underset{n\in 
\mathbb{N}
^{\ast }}{\sup }\left( \left\vert b_{n}\right\vert \right) <+\infty%
\end{array}%
\right.
\end{equation*}%
These conditions rewrite 
\begin{equation*}
\alpha _{n}\in 
\mathbb{R}
^{+},\text{ }n\in 
\mathbb{N}
^{\ast }
\end{equation*}

Now let $\left( \alpha _{n}\right) _{n\in 
\mathbb{N}
^{\ast }}\in l_{\infty }\left( \mathbb{H}\right) $ be a positive element of $%
l_{\infty }\left( \mathbb{H}\right) \mathfrak{.}$ $\left( \alpha _{n}\right)
_{n\in 
\mathbb{N}
^{\ast }}$ is invertible in $l_{\infty }\left( \mathbb{H}\right) $ if and
only if there exists a sequence $\left( c_{n}\right) _{n\in 
\mathbb{N}
^{\ast }}\in l_{\infty }\left( \mathbb{H}\right) $ such that 
\begin{equation*}
\left( \alpha _{n}\right) _{n\in 
\mathbb{N}
^{\ast }}\left( c_{n}\right) _{n\in 
\mathbb{N}
^{\ast }}=\left( 1\right) _{n\in 
\mathbb{N}
^{\ast }}
\end{equation*}%
That is 
\begin{equation*}
\left\{ 
\begin{array}{c}
\alpha _{n}c_{n}=1,\text{ }n\in 
\mathbb{N}
^{\ast } \\ 
\underset{n\in 
\mathbb{N}
^{\ast }}{\sup }\left( \left\vert c_{n}\right\vert \right) <+\infty%
\end{array}%
\right.
\end{equation*}%
These conditions rewrite 
\begin{equation*}
\underset{n\in 
\mathbb{N}
^{\ast }}{\sup }\left( \frac{1}{\alpha _{n}}\right) <+\infty
\end{equation*}%
That is 
\begin{equation*}
0<\underset{n\in 
\mathbb{N}
^{\ast }}{\inf }\left( \alpha _{n}\right)
\end{equation*}

The conclusion is that $\left( \alpha _{n}\right) _{n\in 
\mathbb{N}
^{\ast }}$ is a strictly positive element of $l_{\infty }\left( \mathbb{H}%
\right) $ if and only if%
\begin{equation*}
\left\{ 
\begin{array}{c}
\alpha _{n}\in 
\mathbb{R}
^{+\ast },\text{ }n\in 
\mathbb{N}
^{\ast } \\ 
0<\underset{n\in 
\mathbb{N}
^{\ast }}{\inf }\left( \alpha _{n}\right)%
\end{array}%
\right.
\end{equation*}

The final conclusion is that a sequence $\left( \left( a_{n,k}\right) _{k\in 
\mathbb{N}
^{\ast }}\right) _{n\in 
\mathbb{N}
^{\ast }}$of $\left( l_{\infty }\left( \mathbb{H}\right) \right) ^{%
\mathbb{N}
^{\ast }}$ belongs to $\mathcal{W}\left( l_{\infty }\left( \mathbb{H}\right)
\right) $ if and only if%
\begin{equation*}
\left\{ 
\begin{array}{c}
a_{n,k}\in 
\mathbb{R}
^{+\ast },\text{ }n,\text{ }k\in 
\mathbb{N}
^{\ast } \\ 
0<\underset{n\in 
\mathbb{N}
^{\ast }}{\inf }\left( a_{n,k}\right) ,\text{ }k\in 
\mathbb{N}
^{\ast }%
\end{array}%
\right.
\end{equation*}

The proof of the proposition is then achieved.

$\blacksquare $

\bigskip

We denote by $\preccurlyeq ,$ as usual the partial order on the C*-algebra $%
l_{\infty }\left( \mathbb{H}\right) $ defined by positive elements of $%
l_{\infty }\left( \mathbb{H}\right) \mathfrak{.}$ In view of the last
proposition, it is then clear that for every $a:=\left( a_{n}\right) _{n\in 
\mathbb{N}
^{\ast }},$ $b:=\left( b_{n}\right) _{n\in 
\mathbb{N}
^{\ast }}\in l_{\infty }\left( \mathbb{H}\right) \mathfrak{,}$ we have $%
a\preccurlyeq b$ if and only if $b_{n}-a_{n}\in 
\mathbb{R}
^{+}$ for every $n\in 
\mathbb{N}
^{\ast }$.

\bigskip

\textbf{Proposition 4.6.}

\textit{A sequence }$\left( a_{n}\right) _{n\in 
\mathbb{N}
^{\ast }}$ \textit{of elements of} $\mathcal{W}\left( l_{\infty }\left( 
\mathbb{H}\right) \right) $ \textit{belongs to} $\mathfrak{m}\left( \left(
U_{n}\right) _{n\in 
\mathbb{N}
^{\ast }}\right) $ \textit{if and only if }%
\begin{equation*}
\left\{ 
\begin{array}{c}
0<\underset{n=1}{\overset{+\infty }{\sum }}a_{n,k}^{2}d_{k,n},\text{ }k\in 
\mathbb{N}
^{\ast },\text{ }k\in 
\mathbb{N}
^{\ast } \\ 
\underset{k=1}{\overset{+\infty }{\sum }}\underset{n=1}{\overset{+\infty }{%
\sum }}a_{n,k}^{2}d_{k,n}<+\infty%
\end{array}%
\right.
\end{equation*}

\textit{In this case }$\left( \left( U_{n},a_{n}\right) _{n\in 
\mathbb{N}
^{\ast }}\right) $\textit{\ is a } $\left( \sqrt{\underset{n=1}{\overset{%
+\infty }{\sum }}a_{n,k}^{2}d_{k,n}}\right) _{k\in 
\mathbb{N}
^{\ast }}$\textit{-tight *-fusion frame of }$\mathfrak{\left( \mathbb{H}%
\right) }_{l\mathfrak{_{2}}}\mathfrak{.}$

\bigskip

\textbf{Proof}

Let $\left( a_{n,k}\right) _{n\in 
\mathbb{N}
^{\ast }}$ $\in \mathcal{W}\left( l_{\infty }\left( \mathbb{H}\right)
\right) .$

$\left( a_{n}\right) _{n\in 
\mathbb{N}
^{\ast }}:=\left( \left( a_{n,k}\right) _{k\in 
\mathbb{N}
^{\ast }}\right) _{n\in 
\mathbb{N}
^{\ast }}\in $ $\mathfrak{m}\left( \left( U_{n}\right) _{n\in 
\mathbb{N}
^{\ast }}\right) $ if and only if there exists two elements $\alpha :=\left(
\alpha _{n}\right) _{n\in 
\mathbb{N}
^{\ast }}$ and $\beta :=\left( \beta _{n}\right) _{n\in 
\mathbb{N}
^{\ast }}$ belonging to $\mathcal{W}\left( l_{\infty }\left( \mathbb{H}%
\right) \right) $ such that

i. for every $x:=\left( x_{n}\right) _{n\in 
\mathbb{N}
^{\ast }}$ $\in \left( \mathbb{H}\right) _{l_{2}}\mathfrak{,}$ the series $%
\sum a_{n}^{2}\left\langle P_{U_{n}}\left( x\right) ,P_{U_{n}}\left(
x\right) \right\rangle $ is convergent in $l_{\infty }\left( \mathbb{H}%
\right) \mathfrak{;}$

ii. the following relation holds%
\begin{equation*}
\alpha ^{2}\left\langle x,x\right\rangle \preccurlyeq \underset{n=1}{\overset%
{+\infty }{\sum }}a_{n}^{2}\left\langle P_{U_{n}}\left( x\right)
,P_{U_{n}}\left( x\right) \right\rangle \preccurlyeq \beta ^{2}\left\langle
x,x\right\rangle
\end{equation*}%
But%
\begin{equation*}
a_{n}^{2}\left\langle P_{U_{n}}\left( x\right) ,P_{U_{n}}\left( x\right)
\right\rangle =\left( a_{n,k}^{2}d_{n,k}\left\vert x_{k}\right\vert
^{2}\right) _{k\in 
\mathbb{N}
^{\ast }},\text{ }n\in 
\mathbb{N}
^{\ast }
\end{equation*}%
Hence $\left( a_{n}\right) _{n\in 
\mathbb{N}
^{\ast }}\in $ $\mathfrak{m}\left( \left( U_{n}\right) _{n\in 
\mathbb{N}
^{\ast }}\right) $ if and only if

i'. the series $\sum \left( a_{n,k}^{2}d_{n,k}\left\vert x_{k}\right\vert
^{2}\right) _{k\in 
\mathbb{N}
^{\ast }}$ is convergent in $\left( \mathbb{H}\right) _{l_{2}}$ for every $%
x:=\left( x_{k}\right) _{k\in 
\mathbb{N}
^{\ast }}\in \left( \mathbb{H}\right) _{l_{2}}\mathfrak{;}$

ii'. $\left( \alpha _{k}^{2}\left\vert x_{k}\right\vert ^{2}\right) _{k\in 
\mathbb{N}
^{\ast }}\preccurlyeq \underset{n=1}{\overset{+\infty }{\sum }}\left(
a_{n,k}^{2}d_{k,n}\left\vert x_{k}\right\vert ^{2}\right) _{k\in 
\mathbb{N}
^{\ast }}\preccurlyeq \left( \beta _{k}^{2}\left\vert x_{k}\right\vert
^{2}\right) _{k\in 
\mathbb{N}
^{\ast }}.$

It follows then clearly that $\left( a_{n}\right) _{n\in 
\mathbb{N}
^{\ast }}\in $ $\mathfrak{m}\left( \left( U_{n}\right) _{n\in 
\mathbb{N}
^{\ast }}\right) $ if and only if

\begin{equation*}
\left\{ 
\begin{array}{c}
\underset{n=1}{\overset{+\infty }{\sum }}a_{n,k}^{2}d_{n,k}<+\infty ,k\in 
\mathbb{N}
^{\ast } \\ 
\underset{k=1}{\overset{+\infty }{\sum }}\left( \left( \underset{n=1}{%
\overset{+\infty }{\sum }}a_{n,k}^{2}d_{n,k}\right) \left\vert
x_{k}\right\vert ^{2}\right) ^{2}<+\infty \\ 
\alpha _{k}^{2}\leq \underset{n=1}{\overset{+\infty }{\sum }}%
a_{n,k}^{2}d_{k,n}\leq \beta _{k}^{2}\text{, }k\in 
\mathbb{N}
^{\ast }%
\end{array}%
\right.
\end{equation*}%
holds for every $x:=\left( x_{n}\right) _{n\in 
\mathbb{N}
^{\ast }}$ $\in l_{2}\left( \mathbb{H}\right) \mathfrak{.}$ But since the
mapping%
\begin{equation*}
\begin{array}{ccc}
l_{2}\left( 
\mathbb{C}
\right) & \rightarrow & l_{1}\left( 
\mathbb{C}
\right) \\ 
\left( t_{n}\right) _{n\in 
\mathbb{N}
^{\ast }} & \mapsto & \left( t_{n}^{2}\right) _{n\in 
\mathbb{N}
^{\ast }}%
\end{array}%
\end{equation*}%
is onto, we obtain the following equivalence :

$\left( a_{n}\right) _{n\in 
\mathbb{N}
^{\ast }}\in $ $\mathfrak{m}\left( \left( U_{n}\right) _{n\in 
\mathbb{N}
^{\ast }}\right) $ if and only if

\begin{equation*}
\left\{ 
\begin{array}{c}
\underset{n=1}{\overset{+\infty }{\sum }}a_{n,k}^{2}d_{n,k}<+\infty ,\text{ }%
k\in 
\mathbb{N}
^{\ast } \\ 
\underset{k=1}{\overset{+\infty }{\sum }}\left( \left( \underset{n=1}{%
\overset{+\infty }{\sum }}a_{n,k}^{2}d_{n,k}\right) \left\vert
s_{k}\right\vert ^{2}\right) ^{2}<+\infty \\ 
\alpha _{k}^{2}\leq \underset{n=1}{\overset{+\infty }{\sum }}%
a_{n,k}^{2}d_{k,n}\leq \beta _{k}^{2}\text{, }k\in 
\mathbb{N}
^{\ast }%
\end{array}%
\right.
\end{equation*}%
holds for every $s:=\left( s_{n}\right) _{n\in 
\mathbb{N}
^{\ast }}$ $\in l_{1}\left( 
\mathbb{C}
\right) .$

Finally it follows easily that $\left( a_{n}\right) _{n\in 
\mathbb{N}
^{\ast }}\in $ $\mathfrak{m}\left( \left( U_{n}\right) _{n\in 
\mathbb{N}
^{\ast }}\right) $ if and only if 
\begin{equation*}
\left\{ 
\begin{array}{c}
0<\underset{n=1}{\overset{+\infty }{\sum }}a_{n,k}^{2}d_{n,k}<+\infty ,\text{
}k\in 
\mathbb{N}
^{\ast } \\ 
\underset{k\in 
\mathbb{N}
^{\ast }}{\sup }\left( \underset{n=1}{\overset{+\infty }{\sum }}%
a_{n,k}^{2}d_{n,k}\right) <+\infty \\ 
0<\underset{n=1}{\overset{+\infty }{\sum }}a_{n,k}^{2}d_{n,k},\text{ }k\in 
\mathbb{N}
^{\ast } \\ 
\underset{k=1}{\overset{+\infty }{\sum }}\underset{n=1}{\overset{+\infty }{%
\sum }}a_{n,k}^{2}d_{k,n}<+\infty%
\end{array}%
\right.
\end{equation*}%
That is%
\begin{equation*}
\left\{ 
\begin{array}{c}
0<\underset{n=1}{\overset{+\infty }{\sum }}a_{n,k}^{2}d_{k,n},\text{ }k\in 
\mathbb{N}
^{\ast },\text{ }k\in 
\mathbb{N}
^{\ast } \\ 
\underset{k=1}{\overset{+\infty }{\sum }}\underset{n=1}{\overset{+\infty }{%
\sum }}a_{n,k}^{2}d_{k,n}<+\infty%
\end{array}%
\right.
\end{equation*}%
In this case we can take for $\alpha :=\left( \alpha _{n}\right) _{n\in 
\mathbb{N}
^{\ast }}$ and $\beta :=\left( \beta _{n}\right) _{n\in 
\mathbb{N}
^{\ast }}$ the values

\begin{equation*}
\alpha =\beta =\left( \sqrt{\underset{n=1}{\overset{+\infty }{\sum }}%
a_{n,k}^{2}d_{k,n}}\right) _{k\in 
\mathbb{N}
^{\ast }}
\end{equation*}%
and $\left( \left( U_{n},a_{n}\right) \right) _{n\in 
\mathbb{N}
^{\ast }}$ becomes a $\left( \sqrt{\underset{n=1}{\overset{+\infty }{\sum }}%
a_{n,k}^{2}d_{k,n}}\right) _{k\in 
\mathbb{N}
^{\ast }}$-tight *-fusion frame of $\left( \mathbb{H}\right) _{l_{2}}%
\mathfrak{.}$

The proof of the proposition is then complete.

$\blacksquare $

\bigskip

\section{\textbf{Perturbation of a *-fusion frame}}

\bigskip

\subsection{A distance in the set of all orthogonally complemented subspaces
of a Hilbert $\mathfrak{A}$-module $\mathfrak{H}$}

\bigskip

We denote by $\mathfrak{Comp}\left( \mathfrak{H}\right) $ the set of all
orthogonally complemented subspaces of a Hilbert $\mathfrak{A}$-module $%
\mathfrak{H}$. Let $\mathfrak{d}$ be the mapping 
\begin{equation*}
\begin{array}{cccc}
\mathfrak{d:} & \mathfrak{Comp}\left( \mathfrak{H}\right) ^{2} & \rightarrow
& 
\mathbb{R}
\\ 
& \left( U,V\right) & \mapsto & \mathfrak{d}\left( U,V\right) :=\left\Vert
\pi _{U}-\pi _{V}\right\Vert _{End_{\mathfrak{A}}^{\ast }\left( \mathfrak{H}%
\right) }%
\end{array}%
\end{equation*}%
which was already introduced by D. S. Djordjevic in (\cite{DJOR}). It is
clear that $\mathfrak{d}$ is a well-defined mapping and that%
\begin{equation*}
\left\{ 
\begin{array}{c}
\mathfrak{d}\left( U,V\right) \in 
\mathbb{R}
^{+},\text{ }U,\text{ }V\in \mathfrak{Comp}\left( \mathfrak{H}\right) \\ 
\ \mathfrak{d}\left( U,V\right) =\mathfrak{d}\left( V,U\right) ,\text{ }U,%
\text{ }V\in \mathfrak{Comp}\left( \mathfrak{H}\right) \\ 
\mathfrak{d}\left( U,W\right) \leq \mathfrak{d}\left( U,V\right) +\mathfrak{d%
}\left( V,W\right) ,\text{ }U,\text{ }V,\text{ }W\in \mathfrak{Comp}\left( 
\mathfrak{H}\right) .%
\end{array}%
\right.
\end{equation*}%
Furthermore we have for every $U,$ $V\in \mathfrak{Comp}\left( \mathcal{H}%
\right) $ 
\begin{eqnarray*}
\mathfrak{d}\left( U,V\right) &=&0\implies \pi _{U}=\pi _{V} \\
&\implies &\pi _{U}\left( \mathfrak{H}\right) =\pi _{V}\left( \mathfrak{H}%
\right) \\
&\implies &U=V
\end{eqnarray*}%
The conclusion is that $\mathfrak{d}$ is a distance on the set $\mathfrak{%
Comp}\left( \mathfrak{H}\right) .$

\bigskip

\subsection{The angle between two orthogonally complemented subspaces of a
Hilbert $\mathfrak{A}$-module $\mathfrak{H}$.}

\bigskip

In the mentionned paper (\cite{DJOR}), the author proved that 
\begin{equation}
\mathfrak{d}\left( U,V\right) =\left\Vert \pi _{U}-\pi _{V}\right\Vert
_{End_{\mathfrak{A}}^{\ast }\left( \mathfrak{H}\right) }\leq 1,\text{ }U,%
\text{ }V\in \mathfrak{Comp}\left( \mathfrak{H}\right)  \label{angle}
\end{equation}%
We can then set for every $U,$ $V\in \mathfrak{Comp}\left( \mathfrak{H}%
\right) $ 
\begin{equation*}
\widehat{\left( U,V\right) }:=\arcsin \left( \mathfrak{d}\left( U,V\right)
\right)
\end{equation*}%
The number $\widehat{\left( U,V\right) }$ which belongs to the intervall $%
\left[ 0;\frac{\pi }{2}\right] $ is called the angle between the
orthogonally complemented subspaces $U$ and $V.$

\bigskip

\textbf{Remark}

\textit{Let }$U,$ $V\in \mathfrak{Comp}\left( \mathfrak{H}\right) $ \textit{%
be orthogonal}, \textit{then} $\widehat{\left( U,V\right) }=\frac{\pi }{2}$ 
\textit{but the converse is in general false.}

\bigskip

\textbf{Proof}

1. Assume that $U,$ $V\in \mathfrak{Comp}\left( \mathfrak{H}\right) $ are
orthogonal. Then 
\begin{equation*}
\left\vert \pi _{U}\left( x\right) -\pi _{V}\left( x\right) \right\vert
^{2}=\left\vert \pi _{U}\left( x\right) \right\vert ^{2}+\left\vert \pi
_{V}\left( x\right) \right\vert ^{2},\text{ }x\in \mathfrak{H}
\end{equation*}%
So we have for each $x\in U$%
\begin{eqnarray*}
\left\vert \pi _{U}\left( x\right) -\pi _{V}\left( x\right) \right\vert ^{2}
&=&\left\vert \pi _{U}\left( x\right) \right\vert ^{2} \\
&=&\left\vert x\right\vert ^{2}
\end{eqnarray*}%
It follows that 
\begin{equation*}
\left\Vert \left( \pi _{U}-\pi _{V}\right) \left( x\right) \right\Vert _{%
\mathfrak{H}}=\left\Vert x\right\Vert _{\mathfrak{H}},\text{ }x\in U
\end{equation*}%
Hence we have%
\begin{equation*}
1\leq \left\Vert \pi _{U}-\pi _{V}\right\Vert _{End_{\mathfrak{A}}^{\ast
}\left( \mathfrak{H}\right) }
\end{equation*}%
Consequently, in view of (\ref{angle}), we obtain 
\begin{equation*}
\left\Vert \pi _{U}-\pi _{V}\right\Vert _{End_{\mathfrak{A}}^{\ast }\left( 
\mathfrak{H}\right) }=1
\end{equation*}%
So 
\begin{equation*}
\widehat{\left( U,V\right) }=\frac{\pi }{2}
\end{equation*}

2. Let us consider the classical Hilbert space $l_{2}\left( 
\mathbb{C}
\right) $ which is a Hilbert $%
\mathbb{C}
$-module, $%
\mathbb{C}
$ being viewed as a C*-algebra. Let us set%
\begin{equation*}
\left\{ 
\begin{array}{c}
U_{0}:=\left\{ \left( x_{n}\right) _{n\in 
\mathbb{N}
^{\ast }}\in l_{2}\left( 
\mathbb{C}
\right) :x_{n}=0\text{ if }n\in 
\mathbb{N}
^{\ast }\backslash 4%
\mathbb{N}
\right\} \\ 
V_{0}:=\left\{ \left( x_{n}\right) _{n\in 
\mathbb{N}
^{\ast }}\in l_{2}\left( 
\mathbb{C}
\right) :x_{n}=0\text{ if }n\in 
\mathbb{N}
^{\ast }\backslash 2%
\mathbb{N}
\right\}%
\end{array}%
\right.
\end{equation*}%
It is clear that $U_{0}$ and $V_{0}$ are an orthogonally Hilbert $%
\mathbb{C}
$-submodules of the Hilbert $%
\mathbb{C}
$-module $l_{2}\left( 
\mathbb{C}
\right) $ which are not orthogonal and that 
\begin{equation*}
\left( P_{U_{0}}-P_{V_{0}}\right) \left( \left( 0,1,0,0,....\right) \right)
=\left( 0,1,0,0,....\right)
\end{equation*}%
Hence 
\begin{equation*}
\left\Vert P_{U_{0}}-P_{V_{0}}\right\Vert _{End_{\mathfrak{%
\mathbb{C}
}}^{\ast }\left( l_{2}\left( 
\mathbb{C}
\right) \right) }=1
\end{equation*}%
It follows that 
\begin{equation*}
\widehat{\left( U_{0},V_{0}\right) }=\frac{\pi }{2}
\end{equation*}%
Consequently $U_{0},V_{0}\in \mathfrak{Comp}\left( l_{2}\left( 
\mathbb{C}
\right) \right) $ are not orthogonal but $\widehat{\left( U_{0},V_{0}\right) 
}=\frac{\pi }{2}.\square $

\bigskip

\subsection{Ecart on the set $\mathfrak{Comp}\left( \mathfrak{H}\right) ^{%
\mathbb{N}
^{\ast }}$ of all\ the sequences of orthogonally complemented subspaces of a
Hilbert $\mathfrak{A}$-module $\mathfrak{H}$}

\bigskip

Let $w:=\left( w_{n}\right) _{n\in 
\mathbb{N}
^{\ast }}$ be a sequence of strictly positive real numbers. We consider the
mapping 
\begin{equation*}
\begin{array}{cccc}
\underset{w}{\mathfrak{d}}: & \left( \mathfrak{Comp}\left( \mathfrak{H}%
\right) ^{%
\mathbb{N}
^{\ast }}\right) ^{2} & \rightarrow & 
\mathbb{R}
\cup \left\{ +\infty \right\} \\ 
& \left( \left( U_{n}\right) _{n\in 
\mathbb{N}
^{\ast }},\left( V_{n}\right) _{n\in 
\mathbb{N}
^{\ast }}\right) & \mapsto & \sqrt{\underset{n=1}{\overset{+\infty }{\sum }}%
w_{n}\mathfrak{d}\left( U_{n},V_{n}\right) ^{2}}%
\end{array}%
\end{equation*}%
The mapping $\underset{w}{\mathfrak{d}}$ is well-defined and is an ecart on
the set $\mathfrak{Comp}\left( \mathfrak{H}\right) ^{%
\mathbb{N}
^{\ast }}$((\cite{CHOQ}), pages 61--64)$.$ Indeed it is easy to prove that $%
\underset{w}{\mathfrak{d}}$ fullfiles, for every $\left( U_{n}\right) _{n\in 
\mathbb{N}
^{\ast }},$ $\left( V_{n}\right) _{n\in 
\mathbb{N}
^{\ast }},$ $\left( W_{n}\right) _{n\in 
\mathbb{N}
^{\ast }}$ $\in \mathfrak{Comp}\left( \mathfrak{H}\right) ^{%
\mathbb{N}
^{\ast }},$ the following properties%
\begin{equation*}
\left\{ 
\begin{array}{c}
\underset{w}{\mathfrak{d}}\left( \left( U_{n}\right) _{n\in 
\mathbb{N}
^{\ast }},\left( V_{n}\right) _{n\in 
\mathbb{N}
^{\ast }}\right) \in 
\mathbb{R}
^{+}\cup \left\{ +\infty \right\} \\ 
\ \underset{w}{\mathfrak{d}}\left( \left( U_{n}\right) _{n\in 
\mathbb{N}
^{\ast }},\left( V_{n}\right) _{n\in 
\mathbb{N}
^{\ast }}\right) =\text{ }\underset{w}{\mathfrak{d}}\left( \left(
V_{n}\right) _{n\in 
\mathbb{N}
^{\ast }},\left( U_{n}\right) _{n\in 
\mathbb{N}
^{\ast }}\right) \\ 
\underset{w}{\mathfrak{d}}\left( \left( U_{n}\right) _{n\in 
\mathbb{N}
^{\ast }},\left( W_{n}\right) _{n\in 
\mathbb{N}
^{\ast }}\right) \leq \text{ }\underset{w}{\mathfrak{d}}\left( \left(
U_{n}\right) _{n\in 
\mathbb{N}
^{\ast }},\left( V_{n}\right) _{n\in 
\mathbb{N}
^{\ast }}\right) +\underset{w}{\mathfrak{d}}\left( \left( V_{n}\right)
_{n\in 
\mathbb{N}
^{\ast }},\left( W_{n}\right) _{n\in 
\mathbb{N}
^{\ast }}\right)%
\end{array}%
\right.
\end{equation*}%
It follows that $\underset{w}{\mathfrak{d}}$ defines on the set $\mathfrak{%
Comp}\left( \mathfrak{H}\right) ^{%
\mathbb{N}
^{\ast }}$ a topology $\underset{w}{T}$ in the same way as for a distance :

i. we define the open ball $B_{\underset{w}{\mathfrak{d}}}(\left(
U_{n}\right) _{n\in 
\mathbb{N}
^{\ast }},r)$ of center $\left( U_{n}\right) _{n\in 
\mathbb{N}
^{\ast }}\in \mathfrak{Comp}\left( \mathfrak{H}\right) ^{%
\mathbb{N}
^{\ast }}$ and radius $r\in 
\mathbb{R}
^{+}$ 
\begin{equation*}
B_{\underset{w}{\mathfrak{d}}}(\left( U_{n}\right) _{n\in 
\mathbb{N}
^{\ast }},r)=\{\left( V_{n}\right) _{n\in 
\mathbb{N}
^{\ast }}\in \mathfrak{Comp}\left( \mathfrak{H}\right) ^{%
\mathbb{N}
^{\ast }}:\underset{w}{\mathfrak{d}}(\left( U_{n}\right) _{n\in 
\mathbb{N}
^{\ast }},\left( V_{n}\right) )<r\}
\end{equation*}

ii. we define the open sets of $\mathfrak{Comp}\left( \mathfrak{H}\right) ^{%
\mathbb{N}
^{\ast }}$ as the unions of arbitrary families of open balls.

We can easily prove that the topology $\underset{w}{T}$ satisfies the
Hausdorff separation axiom. Hence $\left( \mathfrak{Comp}\left( \mathfrak{H}%
\right) ^{%
\mathbb{N}
^{\ast }},\underset{w}{T}\right) $ is a Hausdorff topological space. Let us
observe that if\ $w\in l_{2}\left( 
\mathbb{C}
\right) $ then $\underset{w}{\mathfrak{d}}$ is real-valued, so the
restriction $\left. \underset{w}{\mathfrak{d}}\right\vert _{\mathfrak{Comp}%
\left( \mathfrak{H}\right) ^{%
\mathbb{N}
^{\ast }}}^{%
\mathbb{R}
}$ is a distance on $\mathfrak{Comp}\left( \mathfrak{H}\right) ^{%
\mathbb{N}
^{\ast }}$ and $\underset{w}{T}$ becomes a metric space topology on $%
\mathfrak{Comp}\left( \mathfrak{H}\right) ^{%
\mathbb{N}
^{\ast }}.$

Let $\omega :=\left( \omega _{n}\right) _{n\in 
\mathbb{N}
^{\ast }}$ $\in $\textit{\ }$\mathcal{W}\left( \mathfrak{A}\right) .$ We
denote by $\mathfrak{q}\left( \omega \right) $ the sequence $\left(
\left\Vert \omega _{n}\right\Vert _{\mathfrak{A}}^{2}\right) _{n\in 
\mathbb{N}
^{\ast }}.$ We represent the set of all $\left( \mathfrak{K}_{n}\right)
_{n\in 
\mathbb{N}
^{\ast }}\in \mathfrak{Comp}\left( \mathfrak{H}\right) ^{%
\mathbb{N}
^{\ast }}$ such that $\left( \left( \mathfrak{K}_{n},\omega _{n}\right)
\right) _{n\in 
\mathbb{N}
^{\ast }}$ is a *-fusion frame by the notation $\mathfrak{Comp}_{\omega
}\left( \mathfrak{H}\right) ^{%
\mathbb{N}
^{\ast }}.$

\bigskip

\bigskip

\subsection{Results on the perturbation of *-fusion frames}

\bigskip

\textbf{Theorem 5.1.}

1. \textit{Let }$\omega :=\left( \omega _{n}\right) _{n\in 
\mathbb{N}
^{\ast }}$ $\in $\textit{\ }$\mathcal{W}\left( \mathfrak{A}\right) $ \textit{%
and} $\left( \mathfrak{H}_{n}\right) _{n\in 
\mathbb{N}
^{\ast }},$ $\left( \mathfrak{K}_{n}\right) _{n\in 
\mathbb{N}
^{\ast }}\in \mathfrak{Comp}\left( \mathfrak{H}\right) ^{%
\mathbb{N}
^{\ast }}.$

\textit{If }$\left( \left( \mathfrak{H}_{n},\omega _{n}\right) \right)
_{n\in 
\mathbb{N}
^{\ast }}$ \textit{is a *-fusion frame of }$\mathfrak{H}$ \textit{of lower
bound} $A\in \mathfrak{A,}\ $\textit{and} \textit{the following condition
holds} 
\begin{equation*}
\underset{\mathfrak{q}\left( \omega \right) }{\mathfrak{d}}\left( \left( 
\mathfrak{H}_{n}\right) _{n\in 
\mathbb{N}
^{\ast }},\left( \mathfrak{K}_{n}\right) _{n\in 
\mathbb{N}
^{\ast }}\right) <\left\Vert A^{-1}\right\Vert _{\mathfrak{A}}^{-1}
\end{equation*}%
\textit{then} $\left( \left( \mathfrak{K}_{n},\omega _{n}\right) \right)
_{n\in 
\mathbb{N}
^{\ast }}$ \textit{will be a *-fusion frame of} $\mathfrak{H.}$

2. $\mathfrak{Comp}_{\omega }\left( \mathfrak{H}\right) ^{%
\mathbb{N}
^{\ast }}$ \textit{is an open set for the topological space} $\left( 
\mathfrak{Comp}\left( \mathfrak{H}\right) ^{%
\mathbb{N}
^{\ast }},\underset{\mathfrak{q}\left( \omega \right) }{T}\right) .$

\bigskip

\textbf{Proof}

1. Assume that $\left( \left( \mathfrak{H}_{n},\omega _{n}\right) \right)
_{n\in 
\mathbb{N}
^{\ast }}$ is a *-fusion frame of $\mathfrak{H}$ of lower bound $A\in 
\mathfrak{A}$ and upper bound $B\in \mathfrak{A},\ $and the following
condition holds 
\begin{equation*}
\underset{\mathfrak{q}\left( \omega \right) }{\mathfrak{d}}\left( \left( 
\mathfrak{H}_{n}\right) _{n\in 
\mathbb{N}
^{\ast }},\left( \mathfrak{K}_{n}\right) _{n\in 
\mathbb{N}
^{\ast }}\right) <\left\Vert A^{-1}\right\Vert _{\mathfrak{A}}
\end{equation*}

For each $n,m\in 
\mathbb{N}
^{\ast },$ $x\in \mathfrak{H}$ we have%
\begin{eqnarray*}
&&\underset{j=n}{\overset{n+m}{\sum }}\omega _{j}^{2}\left\vert P_{\mathfrak{%
K}_{j}}\left( x\right) \right\vert ^{2} \\
&=&\underset{j=n}{\overset{n+m}{\sum }}\omega _{j}^{2}\left\vert P_{%
\mathfrak{K}_{j}}\left( x\right) -P_{\mathfrak{H}_{j}}\left( x\right)
\right\vert ^{2}+\underset{j=n}{\overset{n+m}{\sum }}\omega _{j}^{2}\left(
\left\langle P_{\mathfrak{H}_{j}}\left( x\right) ,P_{\mathfrak{K}_{j}}\left(
x\right) -P_{\mathfrak{H}_{j}}\left( x\right) \right\rangle \right) + \\
&&+\underset{j=n}{\overset{n+m}{\sum }}\omega _{j}^{2}\left( \left\langle P_{%
\mathfrak{K}_{j}}\left( x\right) -P_{\mathfrak{H}_{j}}\left( x\right) ,P_{%
\mathfrak{H}_{j}}\left( x\right) \right\rangle \right) +\underset{j=n}{%
\overset{n+m}{\sum }}\omega _{j}^{2}\left\vert P_{\mathfrak{H}_{j}}\left(
x\right) \right\vert ^{2}
\end{eqnarray*}%
But we know, thanks to proposition, that 
\begin{eqnarray*}
&&\max \left( \left\Vert \underset{j=n}{\overset{n+m}{\sum }}\omega
_{j}^{2}\left( \left\langle P_{\mathfrak{H}_{j}}\left( x\right) ,P_{%
\mathfrak{K}_{j}}\left( x\right) -P_{\mathfrak{H}_{j}}\left( x\right)
\right\rangle \right) \right\Vert _{\mathfrak{A}},\left\Vert \underset{j=n}{%
\overset{n+m}{\sum }}\omega _{j}^{2}\left( \left\langle P_{\mathfrak{K}%
_{j}}\left( x\right) -P_{\mathfrak{H}_{j}}\left( x\right) ,P_{\mathfrak{H}%
_{j}}\left( x\right) \right\rangle \right) \right\Vert _{\mathfrak{A}}\right)
\\
&\leq &\left\Vert \underset{j=n}{\overset{n+m}{\sum }}\omega
_{j}^{2}\left\vert P_{\mathfrak{K}_{j}}\left( x\right) -P_{\mathfrak{H}%
_{j}}\left( x\right) \right\vert ^{2}\right\Vert _{\mathfrak{A}}\left\Vert 
\underset{j=n}{\overset{n+m}{\sum }}\omega _{j}^{2}\left\vert P_{\mathfrak{H}%
_{j}}\left( x\right) \right\vert ^{2}\right\Vert _{\mathfrak{A}}
\end{eqnarray*}%
It follows that%
\begin{eqnarray*}
&&\left\Vert \underset{j=n}{\overset{n+m}{\sum }}\omega _{j}^{2}\left\vert
P_{\mathfrak{K}_{j}}\left( x\right) \right\vert ^{2}\right\Vert _{\mathfrak{A%
}} \\
&\leq &\left\Vert \underset{j=n}{\overset{n+m}{\sum }}\omega
_{j}^{2}\left\vert P_{\mathfrak{K}_{j}}\left( x\right) -P_{\mathfrak{H}%
_{j}}\left( x\right) \right\vert ^{2}\right\Vert _{\mathfrak{A}} \\
&&+2\left\Vert \underset{j=n}{\overset{n+m}{\sum }}\omega _{j}^{2}\left\vert
P_{\mathfrak{K}_{j}}\left( x\right) \right\vert ^{2}\right\Vert _{\mathfrak{A%
}}\left\Vert \underset{j=n}{\overset{n+m}{\sum }}\omega _{j}^{2}\left\langle
P_{\mathfrak{H}_{j}}\left( x\right) ,P_{\mathfrak{H}_{j}}\left( x\right)
\right\rangle \right\Vert _{\mathfrak{A}}+ \\
&&+\left\Vert \underset{j=n}{\overset{n+m}{\sum }}\omega _{j}^{2}\left\vert
P_{\mathfrak{H}_{j}}\left( x\right) \right\vert ^{2}\right\Vert
\end{eqnarray*}%
Hence we obtain%
\begin{eqnarray*}
&&\sqrt{\left\Vert \underset{j=n}{\overset{n+m}{\sum }}\omega
_{j}^{2}\left\vert P_{\mathfrak{K}_{j}}\left( x\right) \right\vert
^{2}\right\Vert _{\mathfrak{A}}}-\sqrt{\left\Vert \underset{j=n}{\overset{n+m%
}{\sum }}\omega _{j}^{2}\left\vert P_{\mathfrak{H}_{j}}\left( x\right)
\right\vert ^{2}\right\Vert _{\mathfrak{A}}} \\
&\leq &\sqrt{\left\Vert \underset{j=n}{\overset{n+m}{\sum }}\omega
_{j}^{2}\left\vert P_{\mathfrak{K}_{j}}\left( x\right) -P_{\mathfrak{H}%
_{j}}\left( x\right) \right\vert ^{2}\right\Vert _{\mathfrak{A}}}
\end{eqnarray*}%
So%
\begin{eqnarray*}
&&\sqrt{\left\Vert \underset{j=n}{\overset{n+m}{\sum }}\omega
_{j}^{2}\left\vert P_{\mathfrak{K}_{j}}\left( x\right) \right\vert
^{2}\right\Vert _{\mathfrak{A}}} \\
&\leq &\sqrt{\left\Vert \underset{j=n}{\overset{n+m}{\sum }}\omega
_{j}^{2}\left\vert P_{\mathfrak{H}_{j}}\left( x\right) \right\vert
^{2}\right\Vert _{\mathfrak{A}}}+\sqrt{\left\Vert \underset{j=n}{\overset{n+m%
}{\sum }}\omega _{j}^{2}\left\vert P_{\mathfrak{K}_{j}}\left( x\right) -P_{%
\mathfrak{H}_{j}}\left( x\right) \right\vert ^{2}\right\Vert _{\mathfrak{A}}}
\end{eqnarray*}%
Since the sequences $\left( \mathfrak{K}_{j}\right) _{j\in 
\mathbb{N}
^{\ast }}$ and $\left( \mathfrak{H}_{j}\right) _{j\in 
\mathbb{N}
^{\ast }}$\ play the same role in the previous reasonning, we have also%
\begin{eqnarray*}
&&\sqrt{\left\Vert \underset{j=n}{\overset{n+m}{\sum }}\omega
_{j}^{2}\left\vert P_{\mathfrak{H}_{j}}\left( x\right) \right\vert
^{2}\right\Vert _{\mathfrak{A}}} \\
&\leq &\sqrt{\left\Vert \underset{j=n}{\overset{n+m}{\sum }}\omega
_{j}^{2}\left\vert P_{\mathfrak{K}_{j}}\left( x\right) \right\vert
^{2}\right\Vert _{\mathfrak{A}}}+\sqrt{\left\Vert \underset{j=n}{\overset{n+m%
}{\sum }}\omega _{j}^{2}\left\vert P_{\mathfrak{K}_{j}}\left( x\right) -P_{%
\mathfrak{H}_{j}}\left( x\right) \right\vert ^{2}\right\Vert _{\mathfrak{A}}}
\end{eqnarray*}%
It follows that%
\begin{eqnarray}
&&\left\vert \sqrt{\left\Vert \underset{j=n}{\overset{n+m}{\sum }}\omega
_{j}^{2}\left\vert P_{\mathfrak{H}_{j}}\left( x\right) \right\vert
^{2}\right\Vert _{\mathfrak{A}}}-\sqrt{\left\Vert \underset{j=n}{\overset{n+m%
}{\sum }}\omega _{j}^{2}\left\vert P_{\mathfrak{K}_{j}}\left( x\right)
\right\vert ^{2}\right\Vert _{\mathfrak{A}}}\right\vert  \label{INQ} \\
&\leq &\sqrt{\left\Vert \underset{j=n}{\overset{n+m}{\sum }}\omega
_{j}^{2}\left\vert P_{\mathfrak{K}_{j}}\left( x\right) -P_{\mathfrak{H}%
_{j}}\left( x\right) \right\vert ^{2}\right\Vert _{\mathfrak{A}}}  \notag
\end{eqnarray}%
But we know that the series $\sum \omega _{j}^{2}\left\vert P_{\mathfrak{H}%
_{j}}\left( x\right) \right\vert ^{2}$ and $\sum \omega _{j}^{2}\left\vert
P_{\mathfrak{K}_{j}}\left( x\right) -P_{\mathfrak{H}_{j}}\left( x\right)
\right\vert ^{2}$ are convergent. It follows that 
\begin{equation*}
\left\{ 
\begin{array}{c}
\underset{n\rightarrow +\infty }{\lim }\underset{m\in 
\mathbb{N}
^{\ast }}{\sup }\left\Vert \underset{j=n}{\overset{n+m}{\sum }}\omega
_{j}^{2}\left\vert P_{\mathfrak{H}_{j}}\left( x\right) \right\vert
^{2}\right\Vert _{\mathfrak{A}}=0 \\ 
\underset{n\rightarrow +\infty }{\lim }\underset{m\in 
\mathbb{N}
^{\ast }}{\sup }\left\Vert \underset{j=n}{\overset{n+m}{\sum }}\omega
_{j}^{2}\left\vert P_{\mathfrak{K}_{j}}\left( x\right) -P_{\mathfrak{H}%
_{j}}\left( x\right) \right\vert ^{2}\right\Vert _{\mathfrak{A}}=0%
\end{array}%
\right.
\end{equation*}%
Consequently the series $\sum \omega _{j}^{2}\left\vert P_{\mathfrak{K}%
_{j}}\left( x\right) \right\vert ^{2}$ is convergent in $\mathfrak{A.}$
Furthermore if we take $n=1$ and tend $m$ to infinity then the relation (\ref%
{INQ}) becomes%
\begin{eqnarray*}
&&\sqrt{\left\Vert \underset{j=1}{\overset{+\infty }{\sum }}\omega
_{j}^{2}\left\vert P_{\mathfrak{K}_{j}}\left( x\right) \right\vert
^{2}\right\Vert _{\mathfrak{A}}}-\sqrt{\left\Vert \underset{j=1}{\overset{%
+\infty }{\sum }}\omega _{j}^{2}\left\vert P_{\mathfrak{H}_{j}}\left(
x\right) \right\vert ^{2}\right\Vert _{\mathfrak{A}}} \\
&\leq &\underset{\mathfrak{q}\left( \omega \right) }{\mathfrak{d}}\left(
\left( \mathfrak{H}_{n}\right) _{n\in 
\mathbb{N}
^{\ast }},\left( \mathfrak{K}_{n}\right) _{n\in 
\mathbb{N}
^{\ast }}\right) \left\Vert x\right\Vert _{\mathfrak{H}}
\end{eqnarray*}%
It follows that%
\begin{equation*}
\left\{ 
\begin{array}{c}
\sqrt{\left\Vert \underset{j=1}{\overset{+\infty }{\sum }}\omega
_{j}^{2}\left\vert P_{\mathfrak{H}_{j}}\left( x\right) \right\vert
^{2}\right\Vert _{\mathfrak{A}}}-\underset{\mathfrak{q}\left( \omega \right) 
}{\mathfrak{d}}\left( \left( \mathfrak{H}_{n}\right) _{n\in 
\mathbb{N}
^{\ast }},\left( \mathfrak{K}_{n}\right) _{n\in 
\mathbb{N}
^{\ast }}\right) \left\Vert x\right\Vert _{\mathfrak{H}}\leq \sqrt{%
\left\Vert \underset{j=1}{\overset{+\infty }{\sum }}\omega
_{j}^{2}\left\vert P_{\mathfrak{K}_{j}}\left( x\right) \right\vert
^{2}\right\Vert _{\mathfrak{A}}} \\ 
\sqrt{\left\Vert \underset{j=1}{\overset{+\infty }{\sum }}\omega
_{j}^{2}\left\vert P_{\mathfrak{K}_{j}}\left( x\right) \right\vert
^{2}\right\Vert _{\mathfrak{A}}}\leq \sqrt{\left\Vert \underset{j=1}{\overset%
{+\infty }{\sum }}\omega _{j}^{2}\left\vert P_{\mathfrak{H}_{j}}\left(
x\right) \right\vert ^{2}\right\Vert _{\mathfrak{A}}}+\underset{\mathfrak{q}%
\left( \omega \right) }{\mathfrak{d}}\left( \left( \mathfrak{H}_{n}\right)
_{n\in 
\mathbb{N}
^{\ast }},\left( \mathfrak{K}_{n}\right) _{n\in 
\mathbb{N}
^{\ast }}\right) \left\Vert x\right\Vert _{\mathfrak{H}}%
\end{array}%
\right.
\end{equation*}%
But we have 
\begin{equation*}
\left\{ 
\begin{array}{c}
\sqrt{\left\Vert \underset{j=1}{\overset{+\infty }{\sum }}\omega
_{j}^{2}\left\vert P_{\mathfrak{H}_{j}}\left( x\right) \right\vert
^{2}\right\Vert _{\mathfrak{A}}}\leq \left\Vert \left\vert Bx\right\vert
\right\Vert \leq \left\Vert B\right\Vert _{\mathfrak{A}}\left\Vert
x\right\Vert _{\mathfrak{H}} \\ 
\left\Vert A^{-1}\right\Vert _{\mathfrak{A}}^{-1}\left\Vert x\right\Vert _{%
\mathfrak{H}}\leq \left\Vert \left\vert Ax\right\vert \right\Vert \leq \sqrt{%
\left\Vert \underset{j=1}{\overset{+\infty }{\sum }}\omega
_{j}^{2}\left\vert P_{\mathfrak{H}_{j}}\left( x\right) \right\vert
^{2}\right\Vert _{\mathfrak{A}}}%
\end{array}%
\right.
\end{equation*}%
Hence%
\begin{equation*}
\left\{ 
\begin{array}{c}
\sqrt{\left\Vert \underset{j=1}{\overset{+\infty }{\sum }}\omega
_{j}^{2}\left\vert P_{\mathfrak{K}_{j}}\left( x\right) \right\vert
^{2}\right\Vert _{\mathfrak{A}}}\leq \left\Vert \left\vert Bx\right\vert
\right\Vert \leq \left( \left\Vert B\right\Vert _{\mathfrak{A}}+\underset{%
\mathfrak{q}\left( \omega \right) }{\mathfrak{d}}\left( \left( \mathfrak{H}%
_{n}\right) _{n\in 
\mathbb{N}
^{\ast }},\left( \mathfrak{K}_{n}\right) _{n\in 
\mathbb{N}
^{\ast }}\right) \right) \left\Vert x\right\Vert _{\mathfrak{H}} \\ 
\left( \left\Vert A^{-1}\right\Vert _{\mathfrak{A}}^{-1}-\underset{\mathfrak{%
q}\left( \omega \right) }{\mathfrak{d}}\left( \left( \mathfrak{H}_{n}\right)
_{n\in 
\mathbb{N}
^{\ast }},\left( \mathfrak{K}_{n}\right) _{n\in 
\mathbb{N}
^{\ast }}\right) \right) \left\Vert x\right\Vert _{\mathfrak{H}}\leq
\left\Vert \left\vert Ax\right\vert \right\Vert \leq \sqrt{\left\Vert 
\underset{j=1}{\overset{+\infty }{\sum }}\omega _{j}^{2}\left\vert P_{%
\mathfrak{K}_{j}}\left( x\right) \right\vert ^{2}\right\Vert _{\mathfrak{A}}}%
\end{array}%
\right.
\end{equation*}%
But $C:=\left\Vert A^{-1}\right\Vert _{\mathfrak{A}}^{-1}-\underset{%
\mathfrak{q}\left( \omega \right) }{\mathfrak{d}}\left( \left( \mathfrak{H}%
_{n}\right) _{n\in 
\mathbb{N}
^{\ast }},\left( \mathfrak{K}_{n}\right) _{n\in 
\mathbb{N}
^{\ast }}\right) >0.$ It follows that%
\begin{equation*}
\left\{ 
\begin{array}{c}
\left\Vert \underset{j=1}{\overset{+\infty }{\sum }}\omega
_{j}^{2}\left\vert P_{\mathfrak{K}_{j}}\left( x\right) \right\vert
^{2}\right\Vert _{\mathfrak{A}}\leq \left( \left\Vert B\right\Vert _{%
\mathfrak{A}}+\underset{\mathfrak{q}\left( \omega \right) }{\mathfrak{d}}%
\left( \left( \mathfrak{H}_{n}\right) _{n\in 
\mathbb{N}
^{\ast }},\left( \mathfrak{K}_{n}\right) _{n\in 
\mathbb{N}
^{\ast }}\right) \right) ^{2}\left\Vert x\right\Vert _{\mathfrak{H}}^{2} \\ 
C^{2}\left\Vert x\right\Vert _{\mathfrak{H}}^{2}\leq \left\Vert \underset{j=1%
}{\overset{+\infty }{\sum }}\omega _{j}^{2}\left\vert P_{\mathfrak{K}%
_{j}}\left( x\right) \right\vert ^{2}\right\Vert _{\mathfrak{A}}%
\end{array}%
\right.
\end{equation*}%
According to theorem 3.6., $\left( \left( \mathfrak{K}_{n},\omega
_{n}\right) \right) _{n\in 
\mathbb{N}
^{\ast }}$ is a *-fusion frame of $\mathfrak{H.}$

2. If $\mathfrak{Comp}_{\omega }\left( \mathfrak{H}\right) ^{%
\mathbb{N}
^{\ast }}=\emptyset ,$ then $\mathfrak{Comp}_{\omega }\left( \mathfrak{H}%
\right) ^{%
\mathbb{N}
^{\ast }}$ is an open set in the topological space $\left( \mathfrak{Comp}%
\left( \mathfrak{H}\right) ^{%
\mathbb{N}
^{\ast }},\underset{\mathfrak{q}\left( \omega \right) }{T}\right) .$

Assume that $\mathfrak{Comp}_{\omega }\left( \mathfrak{H}\right) ^{%
\mathbb{N}
^{\ast }}\neq \emptyset $ and that $\left( \mathfrak{H}_{n}\right) _{n\in 
\mathbb{N}
^{\ast }}\in \mathfrak{Comp}_{\omega }\left( \mathfrak{H}\right) ^{%
\mathbb{N}
^{\ast }}.$ Let $A$ be a lower bound of the *-fusion frame $\left( \left( 
\mathfrak{H}_{n},\omega _{n}\right) \right) _{n\in 
\mathbb{N}
^{\ast }}.$ The fact that we have obtained in the first part of the proof
can be expressed by the inclusion 
\begin{equation*}
B_{\underset{\mathfrak{q}\left( \omega \right) }{\mathfrak{d}}}(\left( 
\mathfrak{H}_{n}\right) _{n\in 
\mathbb{N}
^{\ast }},\left\Vert A^{-1}\right\Vert _{\mathfrak{A}}^{-1})\subset 
\mathfrak{Comp}_{\omega }\left( \mathfrak{H}\right) ^{%
\mathbb{N}
^{\ast }}
\end{equation*}%
Hence $\mathfrak{Comp}_{\omega }\left( \mathfrak{H}\right) ^{%
\mathbb{N}
^{\ast }}$ is an open set in the topological space $\left( \mathfrak{Comp}%
\left( \mathfrak{H}\right) ^{%
\mathbb{N}
^{\ast }},\underset{\mathfrak{q}\left( \omega \right) }{T}\right) .$

The proof of the theorem is then complete.

$\blacksquare $

\bigskip

\textbf{Corollary 5.2.}

\textit{Let }$\omega :=\left( \omega _{n}\right) _{n\in 
\mathbb{N}
^{\ast }}$ $\in $\textit{\ }$\mathcal{W}\left( \mathfrak{A}\right) $ \textit{%
and} $\left( \mathfrak{H}_{n}\right) _{n\in 
\mathbb{N}
^{\ast }},$ $\left( \mathfrak{K}_{n}\right) _{n\in 
\mathbb{N}
^{\ast }}\in \mathfrak{Comp}\left( \mathfrak{H}\right) ^{%
\mathbb{N}
^{\ast }}.$

1. \textit{If }$\left( \left( \mathfrak{H}_{n},\omega _{n}\right) \right)
_{n\in 
\mathbb{N}
^{\ast }}$ \textit{is a *-fusion frame of }$\mathfrak{H}$ \textit{of lower
bound} $A\in \mathfrak{A,}\ $\textit{and} \textit{the following condition
holds} 
\begin{equation*}
\underset{n=1}{\overset{+\infty }{\sum }}\left\Vert \omega _{n}\right\Vert _{%
\mathfrak{A}}^{2}\widehat{\left( U_{n},V_{n}\right) }^{2}<\left\Vert
A^{-1}\right\Vert _{\mathfrak{A}}^{-2}
\end{equation*}%
\textit{then} $\left( \left( \mathfrak{K}_{n},\omega _{n}\right) \right)
_{n\in 
\mathbb{N}
^{\ast }}$ \textit{will be a *-fusion frame of} $\mathfrak{H.}$

2. \textit{Assume that }$\left( \left( \mathfrak{H}_{n},\omega _{n}\right)
\right) _{n\in 
\mathbb{N}
^{\ast }}$ \textit{is a *-fusion frame of }$\mathfrak{H}$ \textit{of lower
bound} $A\in \mathfrak{A,}\ $\textit{and} \textit{that} \textit{the sequence 
}$\omega $\textit{\ is bounded. Then if the following condition holds} 
\begin{equation*}
\underset{n=1}{\overset{+\infty }{\sum }}\widehat{\left( \mathfrak{H}_{n},%
\mathfrak{K}_{n}\right) }^{2}<\left( \frac{\left\Vert A^{-1}\right\Vert _{%
\mathfrak{A}}^{-1}}{\left\Vert \omega \right\Vert _{\infty }}\right) ^{2}
\end{equation*}%
\textit{then} $\left( \left( \mathfrak{K}_{n},\omega _{n}\right) \right)
_{n\in 
\mathbb{N}
^{\ast }}$ \textit{will be a *-fusion frame of} $\mathfrak{H.}$

3. \textit{Assume that }$\left( \left( \mathfrak{H}_{n},\omega _{n}\right)
\right) _{n\in 
\mathbb{N}
^{\ast }}$ \textit{is a *-fusion frame of }$\mathfrak{H}$ \textit{of lower
bound} $A\in \mathfrak{A,}\ $\textit{and} \textit{that} \textit{there is a
real constant }$p\in ]1;+\infty \lbrack $\textit{\ such that the sequence }$%
\mathfrak{q}\left( \omega \right) $\textit{\ belongs to }$l_{p}\left( 
\mathbb{C}
\right) $\textit{. Then if the following condition holds} 
\begin{equation*}
\underset{n=1}{\overset{+\infty }{\sum }}\widehat{\left( \mathfrak{H}_{n},%
\mathfrak{K}_{n}\right) }^{\frac{2p}{p-1}}<\left( \frac{\left\Vert
A^{-1}\right\Vert _{\mathfrak{A}}^{-2}}{\left\Vert \mathfrak{q}\left( \omega
\right) \right\Vert _{p}}\right) ^{\frac{p}{p-1}}
\end{equation*}%
\textit{then} $\left( \left( \mathfrak{K}_{n},\omega _{n}\right) \right)
_{n\in 
\mathbb{N}
^{\ast }}$ \textit{will be a *-fusion frame of} $\mathfrak{H.}$

\bigskip

\textbf{Proof}

Assume that $\left( \left( \mathfrak{H}_{n},\omega _{n}\right) \right)
_{n\in 
\mathbb{N}
^{\ast }}$ is a *-fusion frame of $\mathfrak{H}$ of lower bound $A\in 
\mathfrak{A.}$ We have 
\begin{eqnarray*}
&&\underset{\mathfrak{q}\left( \omega \right) }{\mathfrak{d}}\left( \left( 
\mathfrak{H}_{n}\right) _{n\in 
\mathbb{N}
^{\ast }},\left( \mathfrak{K}_{n}\right) _{n\in 
\mathbb{N}
^{\ast }}\right) \\
&=&\sqrt{\underset{n=1}{\overset{+\infty }{\sum }}\omega _{n}^{2}\mathfrak{d}%
\left( \mathfrak{H}_{n},\mathfrak{K}_{n}\right) ^{2}} \\
&=&\sqrt{\underset{n=1}{\overset{+\infty }{\sum }}\omega _{n}^{2}\sin ^{2}%
\widehat{\left( \mathfrak{H}_{n},\mathfrak{K}_{n}\right) }} \\
&\leq &\sqrt{\underset{n=1}{\overset{+\infty }{\sum }}\omega _{n}^{2}%
\widehat{\left( \mathfrak{H}_{n},\mathfrak{K}_{n}\right) }^{2}}
\end{eqnarray*}

1. Assume that 
\begin{equation*}
\underset{n=1}{\overset{+\infty }{\sum }}\omega _{n}^{2}\widehat{\left( 
\mathfrak{H}_{n},\mathfrak{K}_{n}\right) }^{2}<\left\Vert A^{-1}\right\Vert
_{\mathfrak{A}}^{-2}
\end{equation*}%
It follows that 
\begin{equation*}
\underset{\mathfrak{q}\left( \omega \right) }{\mathfrak{d}}\left( \left( 
\mathfrak{H}_{n}\right) _{n\in 
\mathbb{N}
^{\ast }},\left( \mathfrak{K}_{n}\right) _{n\in 
\mathbb{N}
^{\ast }}\right) <\left\Vert A^{-1}\right\Vert _{\mathfrak{A}}^{-1}
\end{equation*}%
Hence, by virtue of theorem, $\left( \left( \mathfrak{K}_{n},\omega
_{n}\right) \right) _{n\in 
\mathbb{N}
^{\ast }}$ is a *-fusion frame of $\mathfrak{H.}$

2. Assume that the sequence $\mathfrak{q}\left( \omega \right) $\ is bounded
and that 
\begin{equation*}
\underset{n=1}{\overset{+\infty }{\sum }}\widehat{\left( \mathfrak{H}_{n},%
\mathfrak{K}_{n}\right) }^{2}<\frac{\left\Vert A^{-1}\right\Vert _{\mathfrak{%
A}}^{-2}}{\left\Vert \omega \right\Vert _{\infty }^{2}}
\end{equation*}%
Hence%
\begin{eqnarray*}
&&\underset{\mathfrak{q}\left( \omega \right) }{\mathfrak{d}}\left( \left( 
\mathfrak{H}_{n}\right) _{n\in 
\mathbb{N}
^{\ast }},\left( \mathfrak{K}_{n}\right) _{n\in 
\mathbb{N}
^{\ast }}\right) \\
&\leq &\sqrt{\underset{n=1}{\overset{+\infty }{\sum }}\omega _{n}^{2}%
\widehat{\left( \mathfrak{H}_{n},\mathfrak{K}_{n}\right) }^{2}} \\
&\leq &\left\Vert \omega \right\Vert _{\infty }\sqrt{\underset{n=1}{\overset{%
+\infty }{\sum }}\widehat{\left( \mathfrak{H}_{n},\mathfrak{K}_{n}\right) }%
^{2}} \\
&<&\left\Vert A^{-1}\right\Vert _{\mathfrak{A}}^{-1}
\end{eqnarray*}%
Consequently $\left( \left( \mathfrak{K}_{n},\omega _{n}\right) \right)
_{n\in 
\mathbb{N}
^{\ast }}$ is a *-fusion frame of $\mathfrak{H.}$

3. Assume that there is a real constant $p\in ]1;+\infty \lbrack $\ such
that the sequence $\mathfrak{q}\left( \omega \right) $\ belongs to $%
l_{p}\left( 
\mathbb{C}
\right) $\textit{\ }and that 
\begin{equation*}
\underset{n=1}{\overset{+\infty }{\sum }}\widehat{\left( \mathfrak{H}_{n},%
\mathfrak{K}_{n}\right) }^{\frac{2p}{p-1}}<\left( \frac{\left\Vert
A^{-1}\right\Vert _{\mathfrak{A}}^{-2}}{\left\Vert \mathfrak{q}\left( \omega
\right) \right\Vert _{p}}\right) ^{\frac{p}{p-1}}
\end{equation*}%
Hence, applying the well-known H\"{o}lder's inequality, we obtain%
\begin{eqnarray*}
&&\underset{\mathfrak{q}\left( \omega \right) }{\mathfrak{d}}\left( \left( 
\mathfrak{H}_{n}\right) _{n\in 
\mathbb{N}
^{\ast }},\left( \mathfrak{K}_{n}\right) _{n\in 
\mathbb{N}
^{\ast }}\right) \\
&\leq &\sqrt{\underset{n=1}{\overset{+\infty }{\sum }}\omega _{n}^{2}%
\widehat{\left( \mathfrak{H}_{n},\mathfrak{K}_{n}\right) }^{2}} \\
&\leq &\sqrt{\left( \underset{n=1}{\overset{+\infty }{\sum }}\omega
_{n}^{2p}\right) ^{\frac{1}{p}}\left( \underset{n=1}{\overset{+\infty }{\sum 
}}\widehat{\left( \mathfrak{H}_{n},\mathfrak{K}_{n}\right) }^{\frac{2p}{p-1}%
}\right) ^{\frac{p-1}{p}}} \\
&\leq &\sqrt{\text{ }\left\Vert \mathfrak{q}\left( \omega \right)
\right\Vert _{p}\left( \underset{n=1}{\overset{+\infty }{\sum }}\widehat{%
\left( \mathfrak{H}_{n},\mathfrak{K}_{n}\right) }^{\frac{2p}{p-1}}\right) ^{%
\frac{p-1}{p}}} \\
&<&\left\Vert A^{-1}\right\Vert _{\mathfrak{A}}^{-1}
\end{eqnarray*}%
Consequently $\left( \left( \mathfrak{K}_{n},\omega _{n}\right) \right)
_{n\in 
\mathbb{N}
^{\ast }}$ is a *-fusion frame of $\mathfrak{H.}$

$\blacksquare $

\bigskip

Department of Mathematics, Ibn Tofail University, BP 242, Kenitra, Morocco.
Laboratoire \'{E}quations aux D\'{e}riv\'{e}es Partielles, Alg\`{e}bre et G%
\'{e}om\'{e}trie Spectrales.

E-mail address : nadia.assila@uit.ac.ma

\bigskip

Department of Mathematics, Ibn Tofail University, BP 242, Kenitra, Morocco.
Laboratoire \'{E}quations aux D\'{e}riv\'{e}es Partielles, Alg\`{e}bre et G%
\'{e}om\'{e}trie Spectrales.

E-mail address : samkabbaj@yahoo.fr

\bigskip

Department of Mathematics, Ibn Tofail University, BP 242, Kenitra, Morocco.
Laboratoire \'{E}quations aux D\'{e}riv\'{e}es Partielles, Alg\`{e}bre et G%
\'{e}om\'{e}trie Spectrales.

E-mail address : hzoubeir2014@gmail.com

\bigskip

\bigskip

\bigskip

\bigskip

\bigskip

\bigskip

\bigskip

\end{document}